\documentclass[]{article}

\usepackage{amsmath, amsfonts, amsthm, amssymb, latexsym, relsize, float, caption, subcaption, longtable, hhline, siunitx, verbatim, accents, tabularx, mathtools}
\usepackage[titletoc, title]{appendix}
\usepackage[dvipsnames]{xcolor}
\usepackage{enumitem}
\usepackage{multirow}
\usepackage[normalem]{ulem}

\usepackage{graphicx}
\graphicspath{{./figures/}}

\RequirePackage{algorithm}[0.1] 
\usepackage{algorithmic} 

\usepackage{tikz}
\usetikzlibrary{decorations.pathreplacing}
\usetikzlibrary{arrows.meta}

\newcommand{\N}{{\mathbb{N}}}  
\newcommand{\R}{{\mathbb{R}}}  

\newcommand{\bm}[1]{\text{\boldmath $#1$\unboldmath}}
\newcommand{\mat}[1]{\mathbf{#1}}
\newcommand{\bx}{\bm{x}}
\newcommand{\bb}{\bm{b}}
\newcommand{\bn}{\bm{n}}
\newcommand{\balpha}{\bm{\alpha}}
\newcommand{\bmu}{\bm{\mu}}
\newcommand{\bLambda}{\bm{\Lambda}}
\newcommand{\upgd}{u^{{\texttt{PGD}}}}
\newcommand{\vpgd}{v^{{\texttt{PGD}}}}
\newcommand{\deV}{\delta\!v}
\newcommand{\Apgd}{\mathcal{A}^{{\texttt{PGD}}}}
\newcommand{\Fpgd}{\mathcal{F}^{{\texttt{PGD}}}}
\newcommand{\Nfem}{N^{{\texttt{FEM}}}}
\newcommand{\bupgd}{\mathbf{u}^{{\texttt{PGD}}}}
\newcommand{\bvpgd}{\mathbf{v}^{{\texttt{PGD}}}}
\newcommand{\Map}{\bm{\mathcal{M}}_{\mu}}
\newcommand{\Jaco}{\mat{J}_{\mu}}
\newcommand{\detJ}{\operatorname{det}(\mat{J}_{\mu})}
\newcommand{\adjJ}{\operatorname{adj}(\mat{J}_{\mu})}
\newcommand{\adjJt}{\operatorname{adj}(\mat{J}_{\mu}^T)}
\newcommand{\hOmega}{\hat{\Omega}}
\newcommand{\hmu}{\hat{\mu}}
\newcommand{\hx}{\hat{x}}
\newcommand{\hy}{\hat{y}}
\newcommand{\hT}{\hat{T}}
\newcommand{\bhx}{\bm{\hat{x}}}
\newcommand{\tSUPG}{\tau^{{\texttt{SUPG}}}}
\newcommand{\Nmu}{N^{\mu}}
\newcommand{\Nlam}{N^{\Lambda}}
\newcommand{\hX}{h_{x}}
\newcommand{\hMu}{h_{\mu}}
\newcommand{\hLam}{h_{\Lambda}}
\newcommand{\bV}{\mathbf{V}}
\newcommand{\bphi}{\bm{\Phi}}
\newcommand{\bpsi}{\bm{\Psi}}
\newcommand{\bxi}{\bm{\Xi}}
\newcommand{\bups}{\bm{\Upsilon}}
\newcommand{\muMin}{\mu_{\text{min}}}
\newcommand{\muMax}{\mu_{\text{max}}}
\newcommand{\LamMin}{\Lambda_{\text{min}}}
\newcommand{\LamMax}{\Lambda_{\text{max}}}

\usepackage{xifthen}
\newcommand{\pateraBlock}[7]{

    \ifthenelse{#2 = 1}{\draw (0, -0.25*\s) -- (1*\s, -0.25*\s);}{\ifthenelse{#2 = 2}{\draw[dashed] (0, -0.27*\s) -- (1*\s, -0.27*\s);}{}}
    \draw (0, -0.25*\s) -- (0, 0);
    \draw (1*\s, -0.25*\s) -- (1*\s, 0);
    
    \ifthenelse{#3 = 1}{\draw (-0.25*\s, 0) -- (-0.25*\s, 1*\s);}{\ifthenelse{#3 = 2}{\draw[dashed] (-0.27*\s, 0) -- (-0.27*\s, 1*\s);}{}}
    \draw (-0.25*\s, 0) -- (0, 0);
    \draw (-0.25*\s, 1*\s) -- (0, 1*\s);

    \ifthenelse{#4 = 1}{\draw (1.25*\s, 0) -- (1.25*\s, 1*\s);}{\ifthenelse{#4 = 2}{\draw[dashed] (1.27*\s, 0) -- (1.27*\s, 1*\s);}{}}
    \draw (1.25*\s, 0) -- (1*\s, 0);
    \draw (1.25*\s, 1*\s) -- (1*\s, 1*\s);

    \ifthenelse{#5 = 1}{\draw (0, 1.25*\s) -- (1*\s, 1.25*\s);}{\ifthenelse{#5 = 2}{\draw[dashed] (0, 1.27*\s) -- (1*\s, 1.27*\s);}{}}
    \draw (0, 1*\s) -- (0, 1.25*\s);
    \draw (1*\s, 1*\s) -- (1*\s, 1.25*\s);
    
    \draw[red!10, fill = red!10] (0, 0) rectangle (1*\s, 1*\s);
    
    \ifthenelse{#6 = 1}{\node at (0.5*\s, 0.65*\s) {$\mu_{#1}$};}{\ifthenelse{#6 = 2}{\node at (0.5*\s, 0.65*\s) {$\hmu$};}}

    \ifthenelse{#7 = 1}{\node at (0.5*\s, 0.35*\s) {$\Omega^b_{#1}$};}{\ifthenelse{#7 = 2}{\node at (0.5*\s, 0.35*\s) {$\hOmega^b_{#1}$};}}

}

\newcommand{\pateraWing}[4]{

    \ifthenelse{#1 = 1}{\draw[gray!5, fill = gray!5] (0, -0.25*\s) rectangle (1*\s, 0);}{\ifthenelse{#1 = 2}{\draw[blue!5, fill = blue!5] (0, -0.25*\s) rectangle (1*\s, 0);}}{\draw[gray!5, fill = gray!5] (0, -0.20*\s) rectangle (1*\s, 0);} 
    \ifthenelse{#2 = 1}{\draw[gray!5, fill = gray!5] (0, 1*\s) rectangle (1*\s, 1.25*\s);}{\ifthenelse{#2 = 2}{\draw[blue!5, fill = blue!5] (0, 1*\s) rectangle (1*\s, 1.25*\s);}}{\draw[gray!5, fill = gray!5] (0, 1*\s) rectangle (1*\s, 1.20*\s);} 
    \ifthenelse{#3 = 1}{\draw[gray!5, fill = gray!5] (-0.25*\s, 0) rectangle (0, 1*\s);}{\ifthenelse{#3 = 2}{\draw[yellow!5, fill = yellow!5] (-0.25*\s, 0) rectangle (0, 1*\s);}}{\draw[gray!5, fill = gray!5] (-0.20*\s, 0) rectangle (0, 1*\s);} 
    \ifthenelse{#4 = 1}{\draw[gray!5, fill = gray!5] (1*\s, 0) rectangle (1.25*\s, 1*\s);}{\ifthenelse{#4 = 2}{\draw[yellow!5, fill = yellow!5] (1*\s, 0) rectangle (1.25*\s, 1*\s);}}{\draw[gray!5, fill = gray!5] (1*\s, 0) rectangle (1.20*\s, 1*\s);} 

}

\newtheorem{rem}{Remark}
\newcommand{\Neum}[2]{\nu(\boldsymbol{\mu})\nabla{#1}\cdot\bn_{#2}}

\newcommand{\rev}[1]{\textcolor{black}{#1}}

\usepackage{a4wide}

\begin{document}

\begin{center}

\begin{Large}
\textbf{An overlapping domain decomposition method for the solution of parametric elliptic problems via proper generalized decomposition}
\end{Large}

\medskip

Marco Discacciati$^1$, Ben J. Evans$^1$, Matteo Giacomini$^{2,3}$

\medskip

${}^1$ Department of Mathematical Sciences, Loughborough University, Epinal Way, LE11~3TU, Loughborough, United Kingdom. m.discacciati@lboro.ac.uk, b.j.evans@lboro.ac.uk.

${}^2$ Laboratori de C\`alcul Numeric (LaC\`aN), E.T.S. de Ingenier\'ia de Caminos, Canales y Puertos, Universitat Polit\`ecnica de Catalunya, Barcelona, Spain.

${}^3$ Centre Internacional de M\`etodes Num\`erics en Enginyeria (CIMNE), Barcelona, Spain. matteo.giacomini@upc.edu.

\end{center}

%
%
%
%
%
%
%
%

\begin{abstract}
A non-intrusive proper generalized decomposition (PGD) strategy, coupled with an overlapping domain decomposition (DD) method, is proposed to efficiently construct surrogate models of parametric linear elliptic problems.
A parametric multi-domain formulation is presented, with local subproblems featuring arbitrary Dirichlet interface conditions represented through the traces of the finite element functions used for spatial discretization at the subdomain level, with no need for additional auxiliary basis functions.
The linearity of the operator is exploited to devise low-dimensional problems with only few active boundary parameters.
An overlapping Schwarz method is used to glue the local surrogate models, solving a linear system for the nodal values of the parametric solution at the interfaces, without introducing Lagrange multipliers to enforce the continuity in the overlapping region.
The proposed DD-PGD methodology relies on a fully algebraic formulation allowing for real-time computation based on the efficient interpolation of the local surrogate models in the parametric space, with no additional problems to be solved during the execution of the Schwarz algorithm.
Numerical results for parametric diffusion and convection-diffusion problems are presented to showcase the accuracy of the DD-PGD approach, its robustness in different regimes and its superior performance with respect to standard high-fidelity DD methods.
\end{abstract}




\section{Introduction}
\label{sec:introduction}

Model order reduction (MOR) techniques~\cite{Chinesta:2017:ECM} represent established methodologies for the solution of multi-queries problems, e.g., parametric partial differential equations (PDEs), arising from computationally intensive applications such as uncertainty quantification, optimization, data assimilation and real-time control~\cite{Gunzburger-PWG-18}.
While these techniques have achieved full maturity, their employment in the construction of digital twins of large-scale, multi-physics, multi-disciplinary, systems is still limited by the computational cost of the high-fidelity simulations required during the offline phase of the reduced order model (ROM).
In this context, the last decade has witnessed a growing interest towards the combination of domain decomposition (DD) methods \cite{Smith:1996,Quarteroni:1999,Toselli:2005,Dolean-DJN-15} with ROMs, to reduce the number of coupled parameters and/or degrees of freedom of parametric surrogate models, see, e.g., the recent reviews~\cite{Buhr:2021,Klawonn-HKLW-21}. This is particularly critical in the context of multi-physics phenomena~\cite{Hesthaven-DH-23}, such as micro-electro-mechanical systems~\cite{Corigliano-CDM-13} and the Stokes-Darcy problem~\cite{Martini:2015:ACM}.

Stemming from the seminal works by Maday and R\o nquist on the reduced basis element (RBE) method~\cite{Maday:2002:JSC,Maday:2004:SISC,Lovgren:2006:M2AN}, many strategies have been proposed to (i) split the computational domain of a parametric PDE into subdomains, possibly characterized by simpler geometrical shapes, (ii) compute local approximations of the parametric solutions and (iii) efficiently recompose them to obtain a surrogate model of the original problem.
Existing approaches in the literature propose a combination of overlapping and non-overlapping strategies in step (i), different MOR techniques for step (ii), including reduced basis (RB), proper orthogonal decomposition (POD) and proper generalized decomposition (PGD), whereas step (iii) is performed via Lagrange multipliers and Schwarz iterations, just to name a few.

In the context of RB, non-overlapping DD techniques have been mainly adopted. Indeed, the RBE method relies on a non-overlapping DD approach, where the local solutions in the subdomains are constructed using the RB method and glued together using a mortar approach.
To reduce the cost of RBE, in~\cite{Patera-HKP-13,Eftang:2013:IJNME}, the authors leverage the idea of static condensation by expressing the local ROM as a function of the reduced approximation of the interface solution, or port, and a strategy to derive optimal port spaces is proposed in~\cite{Patera-SP-16}.
A variation of the RBE method relying on parametric boundary conditions at the subdomain interfaces is presented in~\cite{Iapichino:2016:CMA} and strategies to couple non-conforming meshes in the RBE context have been studied in~\cite{Antonietti:2016:M2AN} and~\cite{Gervasio-ZMGQ-22}.

Other ROM solutions inspired by classical DD techniques have been proposed in the literature: the RB hybrid method~\cite{Iapichino:2012:CMAME} employs a coarse mesh strategy with a global high-fidelity solution to ensure the continuity of normal fluxes across the interfaces; the hierarchical model reduction~\cite{Perotto-PEV-10} approach uses a non-overlapping DD method to glue local ROMs obtained from separated representations in the longitudinal and transversal directions; a stabilized local POD-ROM is constructed with overlapping and non-overlapping penalizations in~\cite{Baiges-BCI-13}; \cite{Maier:2014:ANM} couples a Dirichlet-Neumann method on conforming meshes with local POD solutions; in~\cite{Barnett:2022:Sandia}, overlapping and non-overlapping Schwarz alternating methods are discussed in the context of POD approximations.
Recently, local surrogate models have also been coupled using optimization-based DD approaches.  
This is the case of~\cite{Iollo-IST-23}, where a constrained optimization procedure is formulated by minimizing the $L^2$-norm of the jump of the local surrogate models at the interfaces between overlapping subdomains, and of~\cite{Rozza-PNTBR-22} that proposes a non-overlapping variational PDE-constrained optimization approach which minimizes the $L^2$-norm of the distance between the local ROMs at the interface, with an appropriate penalization for the fluxes.

An alternative approach to reduce the cost of ROMs stems from rethinking the parametric PDE as a multi-scale problem.
For example, the localized RB multi-scale method~\cite{Ohlberger-OS-15} employs a coarse non-overlapping partition to devise local ROMs, which are later coupled using a discontinuous Galerkin ansatz at the interface.
Following the variational multi-scale rationale, \cite{Veroy-DVRU-23} proposes an additive splitting of the solution into a coarse and a fine scale, and, using the fine-scale information, the coarse local surrogate models -- obtained by solving an oversampling problem with random boundary conditions~\cite{Smetana-BS-18} -- are coupled via a reduced interface basis.

Despite the extensive efforts to couple DD methods and ROMs, existing approaches still present different shortcomings.
On the one hand, non-overlapping strategies are currently limited by the high cost of devising a parametric representation of the global solution at the interface, leading to high-dimensional spaces to be explored using surrogate models.
On the other hand, most of the existing MOR techniques rely on intrusive implementations with respect to the high-fidelity solver employed to compute the snapshots and require the solution of additional, low-dimensional, problems during the online phase to evaluate the surrogate model for a new set of parameters.
This is particularly critical in the context of industrial problems, where commercial and proprietary software is commonly employed for simulations and source codes are not accessible.
To tackle these issues, increasing attention has been recently devoted to overlapping DD strategies~\cite{Barnett:2022:Sandia,Iollo-IST-23} and non-intrusive solutions to build local ROMs with DD techniques, e.g., by relying on purely algebraic formulations~\cite{Hoang:2021:CMAME} or by combining POD with radial basis functions~\cite{Pain-XFHNP-19}, Gaussian process regression~\cite{Pain-XHFMHBANP-19} and autoencoders~\cite{Pain-HWTKSNNMSP-22}.
Although purely data-driven approaches offer appealing solutions for non-intrusive surrogate models, it is well known that they may lack physical interpretability. For this reason, non-intrusive solutions incorporating physical information, e.g., via physics-informed neural networks~\cite{Raissi-RPK-19}, have gained increasing attention in recent years, being also successfully coupled with non-overlapping~\cite{Karniadakis-JKK-20,Karniadakis-JK-21} and overlapping~\cite{Li-LTWL-19,NissenMeyer-MMN-21,Dolean-DHMM-23} DD approaches.

An alternative solution to circumvent the above mentioned issues is represented by PGD~\cite{Chinesta-AMCK-06,Chinesta:2014}.
PGD offers a physics-based \emph{a priori} MOR framework,  with an \emph{offline} phase constructing a rank-one approximation with no prior knowledge of the solution and an \emph{online} phase where efficient evaluations of the surrogate model are performed by simple interpolation in the parametric space.  
This allows a seamless integration of the resulting ROM with any full-order solver, without the need for any extra solution step in the online phase.
Indeed, non-intrusive PGD implementations, paired with software such as SAMCEF, Abaqus, OpenFOAM and MSC-Nastran, have been presented in~\cite{Ladeveze-CNLB-16,Zou-ZCDA-18,Tsiolakis-TGSOH-20,Tsiolakis-TGSOH-22,Cavaliere-CZSLD-22}.
Moreover, a fully algebraic, non-intrusive framework -- the so-called \emph{encapsulated PGD} -- has been recently proposed in~\cite{Diez:2020:ACME}.

In the context of PGD, DD strategies were first introduced in~\cite{Nazeer:2014:CM} for the overlapping case and in~\cite{Huerta:2017:IJNME} for the non-overlapping one.
More precisely, the Arlequin method~\cite{Nazeer:2014:CM} constructs a local PGD solution in each subdomain and exploits Lagrange multipliers, defined as separated functions in the overlapping regions, to couple the surrogate models, thus leading to a global system involving both local unknowns and Lagrange multipliers.
The approach in~\cite{Huerta:2017:IJNME} relies on a non-overlapping Dirichlet-Dirichlet method: during the offline phase, the local surrogate models are computed in each subdomain as a function of a suitable representation of the trace of the unknown at the interface; in the online phase, the interface problem is solved to impose the continuity of fluxes.
Due to the separated representation of the PGD solution, it follows that the resulting interface equation is nonlinear, even when the original problem is a linear PDE,  thus requiring an appropriate iterative scheme, such as the Newton-Raphson method.

In this work,  a DD-PGD computational framework is devised in the context of linear elliptic PDEs to remedy the shortcomings of existing approaches.  
The strategy relies on an overlapping Schwarz algorithm, executed online, to couple the local surrogate models constructed offline in each subdomain. Parametric Dirichlet boundary conditions are employed at the subdomain level and the linearity of the operator is exploited to reduce the dimensionality of each local ROM via superimposition.  This allows to devise a set of low-dimensional problems which can be easily parallelized to enhance the performance of the method. Differently from~\cite{Huerta:2017:IJNME},  the proposed method does not require the definition of any auxiliary basis functions at the interfaces, but it can rely, e.g., on the traces of the finite element functions used for the spatial discretization within each subdomain. 
Moreover,  thanks to an \emph{ad hoc} multi-domain reformulation of the original parametric problem,  the coupling in the online phase occurs only at the interfaces,  instead of across the whole overlapping region as in~\cite{Nazeer:2014:CM}.  The continuity of the solution and of its fluxes is indeed guaranteed by imposing the equality of the traces of the local PGD solutions at the interfaces, without Lagrange multipliers (with separated representations) to glue the local ROMs in the entire overlap. This is practically achieved by formulating the overlapping Schwarz algorithm as a parametric interface linear system, which can be efficiently solved in real time in the online phase by standard matrix-free Krylov methods, while the local surrogate models are evaluated via interpolation in the parametric domain, with no extra solution step.
Finally, the proposed methodology provides a physics-based PGD-ROM,  non-intrusive with respect to the high-fidelity spatial solver and featuring a reduced number of interface parameters in each subproblem solved in the offline phase.

The remainder of this paper is structured as follows. Section~\ref{sec:setting} introduces the parametric elliptic problem, its multi-domain formulation and the main idea of the proposed DD-PGD approach. In Sect.~\ref{sec:offline},  the offline phase of the method is presented,  explaining the rationale for constructing local PGD surrogate models with parametric boundary conditions at the interfaces between subdomains,  while reducing the dimensionality of the resulting problem by exploiting the linearity of the underlying operator. The online phase accounting for the parametric overlapping Schwarz method to solve the linear interface system is described in Sect.~\ref{sec:online}, while in Sect.~\ref{sec:results} numerical tests are presented to assess the accuracy, robustness and efficiency of the proposed methodology. Finally,  Sect.~\ref{sec:Conclusions} summarizes the conclusions of this work and two appendices provide technical details on the encapsulated PGD framework and the implementation of the presented approach.

\section{Problem setting and parametric multi-domain formulation}
\label{sec:setting}

Let $\Omega \subset \R^d$ ($d = 1, 2, 3$) be an open bounded domain with Lipschitz boundary $\partial\Omega = \Gamma^D \cup \Gamma^N$, such that $\Gamma^D \cap \Gamma^N = \emptyset$. Let $\bmu =(\mu^1,\ldots,\mu^P) \in \mathcal{P}$ be a tuple of $P \in \N$ problem parameters with $\mathcal{P} = \mathcal{I}^1 \times \dots \times \mathcal{I}^P \subset \R^P$ and each $\mathcal{I}^p$ compact ($p=1,\ldots,P$).
Consider the linear elliptic parametric operator
\begin{equation*}
  L(u(\bmu);\bmu) = -\nabla\cdot(\nu(\bmu)\nabla u(\bmu)) + \boldsymbol{\alpha}(\bmu)\cdot\nabla u(\bmu) + \gamma(\bmu) u(\bmu) \, ,
\end{equation*}
and the parametric boundary value problem: for all $\bmu \in \mathcal{P}$, find $u(\bmu)$ such that 
\begin{equation}
	\label{eq:globalProb}
	\begin{array}{rcll}
		L(u(\bmu);\bmu) &=& s(\bmu) &\quad \text{in } \Omega,\\
		u(\bmu) &=& g^D(\bmu)& \quad \text{on } \Gamma^D,\\
		\Neum{u(\bmu)}{} &=& g^N(\bmu) & \quad \text{on } \Gamma^N,
	\end{array}
\end{equation}
where $s(\bmu)$ denotes the source term and $g^D(\bmu)$ and $g^N(\bmu)$ are given functions that prescribe Dirichlet and Neumann boundary conditions on $\Gamma^D$ and $\Gamma^N$, respectively, with $\bn$ the unit normal vector to $\Gamma^N$, pointing outwards of the domain. Note that all the above material data, physical quantities and boundary conditions are functions of the parameters $\bmu$. For the sake of readability, the domain and its boundary are assumed to be independent of $\bmu$, although the framework presented in this work can be applied also to geometric parameters, as shown in the numerical example of Sect.~\ref{sec:testRozza}.

To guarantee the well posedness of~\eqref{eq:globalProb}, we assume that, for all $\bmu \in \mathcal{P}$, there exists $\nu_0 >0$ such that $\nu(\bmu) \geq \nu_0$, that $\boldsymbol{\alpha}(\bmu), \gamma(\bmu) \in L^\infty(\Omega)$, and that there exists $\gamma_0 \geq 0$ such that
\begin{equation*}
-\frac{1}{2} \nabla \cdot \boldsymbol{\alpha}(\bmu) + \gamma(\bmu) \geq \gamma_0 \qquad \mbox{a.e. in } \Omega \, ,
\end{equation*}
with
\begin{equation*}
\| \boldsymbol{\alpha}(\bmu) \cdot \bn \|_{L^\infty(\Gamma^N)} < \frac{2}{C_t} \min \left( \frac{\nu_0}{2}, \frac{\nu_0}{2} C_\Omega^{-1} + \gamma_0 \right) \, ,
\end{equation*}
and $C_t,C_\Omega>0$ being the trace and Poincar\'e constants, respectively (see, e.g., \cite{Quarteroni:1994}).

\begin{figure}[bht]
	\centering
	\begin{tikzpicture}
		\fill[gray!40] (3, 0) rectangle (4, 3);
        \node at (3.5, 0.5) {$\Omega_{12}$};
				
        \node[black] at (4.3, 2.3) {$\Gamma_1$};
		
        \node[black] at (2.7, 2.3) {$\Gamma_2$};
		
		\draw[black, thick] (3, 0) rectangle (7, 3);
		\node[black] at (5.5, 1.5) {$\Omega_2$};
		\draw[<->, thick, black,dashed]  (3.1, 1.8) -- (6.9, 1.8);
		
		\draw[black, thick] (0, 0) rectangle (4, 3);
		\node[black] at (1.5, 1.5) {$\Omega_1$};
		\draw[<->, thick, black,dashed]  (0.1, 1.2) -- (3.9, 1.2);
		
	\end{tikzpicture}
\caption{Partition of the domain $\Omega$ into two overlapping subdomains.}\label{fig:partition}
\end{figure}
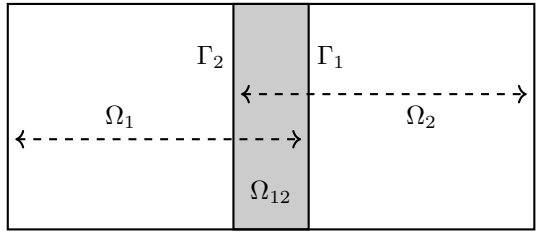

Consider a decomposition of the domain $\Omega$ into two overlapping subdomains $\Omega_i \subset \Omega$ ($i = 1, 2$) such that $\Omega_1 \cup \Omega_2 = \Omega$ and $\Omega_1 \cap \Omega_2 = \Omega_{12} \neq \emptyset$. For $i=1,2$, let $\Gamma_i = \partial\Omega_i \setminus \partial\Omega$ as shown in Fig.~\ref{fig:partition}, let $\Gamma = \cup_{i} \Gamma_i$ be the union of all interfaces, and let $\Gamma_i^D = \Gamma^D \cap \partial\Omega_i$ and $\Gamma_i^N = \Gamma^N \cap \partial\Omega_i$.

For clarity of exposition, we henceforth focus on the case of two subdomains, but the approach can be straightforwardly extended to the case of more than two subdomains without cross-points, as shown in Sect.~\ref{sec:testPatera}.

Problem~\eqref{eq:globalProb} can be rewritten in the equivalent multi-domain formulation:
for all $\bmu \in \mathcal{P}$, find $u_i(\bmu)$ ($i = 1, 2$) such that
\begin{equation}
	\label{eq:multiDom}
	\begin{array}{rcll}
		L(u_i(\bmu);\,\bmu) &=& s_i(\bmu) & \quad \text{in } \Omega_i,\\
		u_i(\bmu) &=& g^D_i(\bmu) & \quad \text{on } \Gamma_i^D,\\
		\Neum{u_i(\bmu)}{} &=& g^N_i(\bmu) & \quad \text{on } \Gamma_i^N,\\
		u_1(\bmu) &=& u_2(\bmu) & \quad \text{on } \Gamma\,,
	\end{array}
\end{equation}
where $s_i(\bmu)$, $g^D_i(\bmu)$ and $g^N_i(\bmu)$ denote the restrictions of $s(\bmu)$, $g^D(\bmu)$ and $g^N(\bmu)$ to $\Omega_i$, $\Gamma_i^D$ and $\Gamma_i^N$, respectively.

It is worth noticing that the continuity of the local solutions $u_1(\bmu)$ and $u_2(\bmu)$ across the interfaces $\Gamma_1$ and $\Gamma_2$ follows from the last condition in~\eqref{eq:multiDom}. 
The equivalence of~\eqref{eq:globalProb} and~\eqref{eq:multiDom} can be proved by straightforwardly extending the argument of Proposition 2.1 of \cite{Discacciati:2013:SICON} to take into account the parameters $\bmu \in \mathcal{P}$. Due to the equivalence, there holds $u_i(\bmu) = u(\bmu)|_{\Omega_i}$ for $i=1,2$ and for all $\bmu \in \mathcal{P}$, and, in particular, $u_1(\bmu) = u_2(\bmu)$ in $\Omega_{12}$.
\rev{Moreover, the continuity of the local solutions in the overlapping region $\Omega_{12}$ straightforwardly guarantees the continuity of the normal fluxes at the interfaces, thus avoiding imposing this condition explicitly}.

\subsection{Parametric overlapping Schwarz method}
The multi-domain problem~\eqref{eq:multiDom} can be used to formulate the overlapping Schwarz method (see, e.g., \cite{Smith:1996}) reported in Algorithm~\ref{alg:Schwarz}. Problems~\eqref{eq:localProblem1} and~\eqref{eq:localProblem2} are solved on the local subdomains $\Omega_1$ and $\Omega_2$, respectively (Algorithm~\ref{alg:Schwarz}, steps 2 and 4), with appropriate Dirichlet-type boundary conditions at the interfaces (Algorithm~\ref{alg:Schwarz}, steps 3 and 5). The algorithm stops when the discrepancy between the local solutions at the interfaces is below a user-defined tolerance $\texttt{tol}$ (Algorithm~\ref{alg:Schwarz}, step 6) and the solution $u(\bar{\bmu})$ of the parametric problem~\eqref{eq:globalProb} for the fixed set of parameters $\bar{\bmu} \in \mathcal{P}$ is retrieved by gluing the local solutions (Algorithm~\ref{alg:Schwarz}, step 7).
For the definition of the stopping criterion, a suitable norm $\|\cdot\|_{\Gamma}$ on $\Gamma$ must be defined. In the numerical simulations of Sect.~\ref{sec:results}, the $l^{\infty}(\Gamma)$ norm is employed.
\begin{subequations}
\begin{algorithm}[!ht]
\caption{Overlapping Schwarz method for a two-domain parametric problem}\label{alg:Schwarz}
\begin{algorithmic}[1]
\REQUIRE{Fixed set of parametric values $\bar{\bmu} \in \mathcal{P}$, tolerance $\texttt{tol}$ for the stopping criterion and initial function $\lambda_1^{(0)}$ defined on the interface $\Gamma_1$.}
\STATE{Set $k=1$.}
\STATE{
Find $u_1^{(k)}\!(\bar{\bmu})$ such that 
	\begin{equation}\label{eq:localProblem1}
		\begin{array}{rcll}
			L (u_1^{(k)}\!(\bar{\bmu}); \, \bar{\bmu}) &=& s_1(\bar{\bmu}) & \quad \text{in } \Omega_1,\\
			u_1^{(k)}\!(\bar{\bmu}) &=& g_1^D(\bar{\bmu}) & \quad \text{on } \Gamma_1^D,\\
			\displaystyle
			\nu(\bar{\bmu}) \nabla u_1^{(k)}\!(\bar{\bmu}) \cdot \bn &=& g_1^N(\bar{\bmu}) & \quad \text{on } \Gamma_1^N,\\
			u_1^{(k)}\!(\bar{\bmu}) &=& \lambda_1^{(k-1)} & \quad \text{on } \Gamma_1 .
		\end{array}
	\end{equation}
	\label{alg:firstStep}
}
\STATE{Set $\lambda_2^{(k)} = u_1^{(k)}\!(\bar{\bmu})|_{\Gamma_2}$.}
\STATE{
Find $u_2^{(k)}\!(\bar{\bmu})$ such that 
	\begin{equation}\label{eq:localProblem2}
		\begin{array}{rcll}
			L (u_2^{(k)}\!(\bar{\bmu}); \, \bar{\bmu}) &=& s_2(\bar{\bmu}) & \quad \text{in } \Omega_2,\\
			u_2^{(k)}\!(\bar{\bmu}) &=& g_2^D(\bar{\bmu}) & \quad \text{on } \Gamma_2^D,\\
			\displaystyle
			\nu(\bar{\bmu}) \nabla u_2^{(k)}\!(\bar{\bmu}) \cdot \bn &=& g_2^N(\bar{\bmu}) & \quad \text{on } \Gamma_2^N,\\
			u_2^{(k)} &=& \lambda_2^{(k)} & \quad \text{on } \Gamma_2.
		\end{array}
	\end{equation}
}
\STATE{Set $\lambda_1^{(k)}= u_2^{(k)}\!(\bar{\bmu})|_{\Gamma_1}$.}
\IF{$\|u_1^{(k)}\!(\bar{\bmu})_{\vert\Gamma} - u_2^{(k)}\!(\bar{\bmu})_{\vert\Gamma}\|_\Gamma < \texttt{tol}$}
\STATE{Construct the solution on the entire domain as
\begin{equation}\label{eq:finalSol}
u(\bar{\bmu}) = \begin{cases}
u_1^{(k)}\!(\bar{\bmu}) & \text{in } \Omega_1 \\
u_2^{(k)}\!(\bar{\bmu}) & \text{in } \Omega_2 \setminus \Omega_{12}
\end{cases}
\end{equation}
}
\ELSE
\STATE{$k \gets k+1$.}
\STATE{\textbf{go to} step~\ref{alg:firstStep}.}
\ENDIF
\ENSURE{Solution $u(\bar{\bmu})$ in the entire domain $\Omega$.}
\end{algorithmic}
\end{algorithm}
\end{subequations}

\subsection{Algebraic formulation of the Schwarz method}
\label{sec:schwarzAlg}

Consider a finite element discretization of problems~\eqref{eq:localProblem1} and~\eqref{eq:localProblem2}.
More precisely, let $\mat{A}_{\Omega_i}$ ($i=1,2$) be the invertible matrix associated with the finite element approximation of the local problem in the subdomain $\Omega_i$,  whose rows and columns correspond to the degrees of freedom inside $\Omega_i$, excluding the unknowns at the interface $\Gamma_i$. 
Moreover, let $\mat{A}_{\Gamma_i}$ ($i=1,2$) be the finite element matrix with rows corresponding to the degrees of freedom inside $\Omega_i$ and columns associated with the degrees of freedom on $\Gamma_i$, whereas let $\mathbf{f}_{\Omega_i}$ ($i=1,2$) denote the finite element vector accounting for the contributions of the source term and the Neumann boundary conditions.
Finally, let $\mat{R}_{\Omega_i\to\Gamma_j}$ ($i,j=1,2$, $i\not= j$) be the restriction matrix that, for any vector of nodal values inside $\Omega_i$, returns the vector of nodal values at the interface $\Gamma_j$ internal to $\Omega_i$. 
Hence, following, e.g., \cite[Sect. 1.1.1]{Smith:1996},  the procedure in Algorithm~\ref{alg:Schwarz} can be rewritten in algebraic form as a block Gauss-Seidel method for the linear system
\begin{equation}\label{eq:systemSchwarz}
\begin{pmatrix}
\mat{A}_{\Omega_1} & \mat{A}_{\Gamma_1} & \mat{0} & \mat{0} \\
\mat{0} & \mat{I}_{\Gamma_1} & - \mat{R}_{\Omega_2\to\Gamma_1} & \mat{0} \\
\mat{0} & \mat{0} & \mat{A}_{\Omega_2} & \mat{A}_{\Gamma_2} \\
- \mat{R}_{\Omega_1\to\Gamma_2} & \mat{0} & \mat{0} & \mat{I}_{\Gamma_2}
\end{pmatrix}
\begin{pmatrix}
\mathbf{u}_{\Omega_1}\!(\bar{\bmu}) \\
\mathbf{u}_{\Gamma_1}\!(\bar{\bmu}) \\
\mathbf{u}_{\Omega_2}\!(\bar{\bmu}) \\
\mathbf{u}_{\Gamma_2}\!(\bar{\bmu})
\end{pmatrix}
=
\begin{pmatrix}
\mathbf{f}_{\Omega_1}\!(\bar{\bmu}) \\
\mathbf{0} \\
\mathbf{f}_{\Omega_2}\!(\bar{\bmu}) \\
\mathbf{0}
\end{pmatrix}\,,
\end{equation}
where $\mathbf{u}_{\Omega_i}(\bar{\bmu})$ and $\mathbf{u}_{\Gamma_i}(\bar{\bmu})$ denote the vectors of nodal values inside the domain $\Omega_i$ and at the interface $\Gamma_i$, respectively, whereas $\mat{I}_{\Gamma_i}$ represents the identity matrix at $\Gamma_i$ ($i=1,2$).

By computing the Schur complement of the linear system~\eqref{eq:systemSchwarz},  the degrees of freedom $\mathbf{u}_{\Omega_1}(\bar{\bmu})$ and $\mathbf{u}_{\Omega_2}(\bar{\bmu})$ internal to each subdomain can be eliminated, expressing them in terms of the unknowns $\mathbf{u}_{\Gamma_1}(\bar{\bmu})$ and $\mathbf{u}_{\Gamma_2}(\bar{\bmu})$ on $\Gamma_1$ and $\Gamma_2$, respectively,  yielding the so-called \emph{interface system}
\begin{equation}\label{eq:systemSchwarzInterface}
\hspace*{-3mm}
\begin{pmatrix}
\mat{I}_{\Gamma_1} & \mat{R}_{\Omega_2\to\Gamma_1} \mat{A}_{\Omega_2}^{-1} \mat{A}_{\Gamma_2} \\[3pt]
\mat{R}_{\Omega_1\to\Gamma_2} \mat{A}_{\Omega_1}^{-1} \mat{A}_{\Gamma_1} & \mat{I}_{\Gamma_2}
\end{pmatrix}
\begin{pmatrix}
\mathbf{u}_{\Gamma_1}\!(\bar{\bmu}) \\[3pt]
\mathbf{u}_{\Gamma_2}\!(\bar{\bmu})
\end{pmatrix}
=
\begin{pmatrix}
\mat{R}_{\Omega_2\to\Gamma_1} \mat{A}_{\Omega_2}^{-1} \mathbf{f}_{\Omega_2}\!(\bar{\bmu}) \\[3pt]
\mat{R}_{\Omega_1\to\Gamma_2} \mat{A}_{\Omega_1}^{-1} \mathbf{f}_{\Omega_1}\!(\bar{\bmu})
\end{pmatrix} \,  ,
\end{equation}
which can be solved using a suitable matrix-free Krylov method, e.g., GMRES~\cite{Saad:1986:SISSC}. 

It is worth noticing that problem~\eqref{eq:systemSchwarzInterface} corresponds to steps 2--5 of Algorithm~\ref{alg:Schwarz}, namely
\begin{equation}\label{eq:SchwarzAlg}
\mat{I}_{\Gamma_j}\mathbf{u}_{\Gamma_j}\!(\bar{\bmu})
=
\mat{R}_{\Omega_i\to\Gamma_j} \mat{A}_{\Omega_i}^{-1} \mathbf{f}_{\Omega_i}\!(\bar{\bmu}) 
+ \mat{R}_{\Omega_i\to\Gamma_j} \mat{A}_{\Omega_i}^{-1} (-\mat{A}_{\Gamma_i}\mathbf{u}_{\Gamma_i}\!(\bar{\bmu})) \, ,
\end{equation}
where the matrix-vector operations on the right-hand side of equation~\eqref{eq:SchwarzAlg} are the algebraic counterpart of the following operations:
\begin{enumerate}[label=(\Alph*)]
\item extension into subdomain $\Omega_i$ of the Dirichlet datum $\mathbf{u}_{\Gamma_i}\!(\bar{\bmu})$ at interface $\Gamma_i$ by the matrix-vector product $-\mat{A}_{\Gamma_i} \mathbf{u}_{\Gamma_i}\!(\bar{\bmu})$;
\item computation of the local solution $\mathbf{u}_{\Omega_i}\!(\bar{\bmu})$ in subdomain $\Omega_i$ as the superposition of $\mat{A}_{\Omega_i}^{-1} \mathbf{f}_{\Omega_i}\!(\bar{\bmu})$ and $\mat{A}_{\Omega_i}^{-1} \left(- \mat{A}_{\Gamma_i}\mathbf{u}_{\Gamma_i}\!(\bar{\bmu}) \right)$, solving a linear system with matrix $\mat{A}_{\Omega_i}$;
\item restriction of the computed solution $\mathbf{u}_{\Omega_i}\!(\bar{\bmu})$ to the internal interface $\Gamma_j$, $j \not= i$, through the restriction matrix $\mat{R}_{\Omega_i\to\Gamma_j}$.
\end{enumerate}

Note that the computational effort required by the above Schwarz method is proportional to the cost of solving the local problems~\eqref{eq:localProblem1} and~\eqref{eq:localProblem2}, that is, the cost of step (B). This can become demanding when a new set of parameters $\bar{\bmu} \in \mathcal{P}$ is to be tested since the entire procedure needs to be executed from scratch. Indeed,  the matrices $\mat{A}_{\Omega_i}$ and $\mat{A}_{\Gamma_i}$ may themselves depend on the parameters $\bar{\bmu}$ defining the novel configuration under analysis.

\subsection{The DD-PGD strategy}
\label{sec:DDPGDstrategy}

To reduce the computational cost, in this paper the Schwarz algorithm is combined with a PGD-based surrogate model to efficiently obtain the solution of the parametric problem~\eqref{eq:globalProb}, for any set of parameters $\bmu \in \mathcal{P}$. To this aim, the Schwarz algorithm is reformulated by identifying an offline phase and an online phase as follows.
\begin{enumerate}
\item
In the \emph{offline} phase, the local parametric problems~\eqref{eq:localProblem1} and~\eqref{eq:localProblem2} are solved using the PGD method to devise a set of surrogate solutions $\upgd_i$ ($i=1,2$) explicitly depending on space, $\bx$, on problem parameters, $\bmu$, and on arbitrary, problem-relevant functions $\lambda_i$\rev{, which represent the traces of the unknown solution on the interfaces $\Gamma_i$.} 
The arbitrariness of the functions $\lambda_i$ is dealt with by parametrizing them through a set of auxiliary parameters, say $\bLambda_i$, as detailed in Sect.~\ref{sec:offline}. The output of the offline phase is the set of local surrogate models $\upgd_i$, featuring arbitrary traces at the interfaces, which are thus suitable for efficient evaluations during the Schwarz algorithm. This procedure is meant to replace the computationally demanding step (B).
\item
In the \emph{online} phase, the Schwarz algorithm is performed using the interface formulation \eqref{eq:systemSchwarzInterface}. For a fixed set of parametric values $\bar{\bmu} \in \mathcal{P}$, at each iteration of the algorithm, the extension of Dirichlet interface data and the solution of the local problems~\eqref{eq:localProblem1} and~\eqref{eq:localProblem2} in steps (A) and (B) are replaced by the evaluation of the precomputed local surrogate models at specific instances of the auxiliary interface parameters $\bLambda_i$. It is worth noticing that, contrary to alternative \emph{a posteriori} ROMs requiring the solution of small problems in the online phase, the evaluation of the PGD surrogate model for a specific value of the parameters only relies on interpolation procedures, thus allowing the Schwarz algorithm to be executed in real time. Details of the online phase are provided in Sect.~\ref{sec:online}.
\end{enumerate}

\rev{
\begin{rem}
The overlapping Schwarz method used in the online phase has been mainly chosen for computational efficiency in the offline phase. Indeed, while non-overlapping DD techniques such as, e.g., Neumann-Neumann or FETI methods~\cite{Toselli:2005} could alternatively be used to rewrite the parametric problem~\eqref{eq:globalProb} into an equivalent multi-domain formulation, in the online phase both the continuity of traces and the continuity of fluxes would have to be imposed through a suitably preconditioned interface equation.
In the offline phase, this would entail the solution of local parametric problems with arbitrary fluxes at the interfaces, besides those with arbitrary traces employed in the present strategy. Therefore, although alternative DD strategies are possible, in this work only the overlapping Schwarz method is considered to avoid increasing the overall computational cost of the DD-ROM procedure.
\end{rem}
}

\section{Local surrogate models using proper generalized decomposition}
\label{sec:offline}

In this section, the procedure to construct the PGD local surrogate models in the offline phase is presented.
The local parametric problem to be solved in the generic subdomain $\Omega_i$ is: for all $\bmu \in \mathcal{P}$, find $u_i(\bmu)$ such that 
\begin{equation}
	\label{eq:localProb}
	\begin{array}{rcll}
		L_{}({u}_{i}(\bmu); \bmu) &=& s_i(\bmu) & \quad \text{in } \Omega_i\,,\\
		u_i(\bmu) &=& g^D_i(\bmu) & \quad \text{on } \Gamma_i^D,\\
		\Neum{u_i(\bmu)}{} &=& g^N_i(\bmu) & \quad \text{on } \Gamma_i^N,\\
		u_i(\bmu) &=& \lambda_i & \quad \text{on } \Gamma_i \, ,
	\end{array}
\end{equation}
for arbitrary Dirichlet data $\lambda_i$ at the interface $\Gamma_i$, where $\lambda_i=\lambda_i(\bx)$ is a space-dependent function. 

The arbitrariness of the boundary function $\lambda_i$ is dealt with by an appropriate parametrization.
Considering that problem \eqref{eq:localProb} is solved in a finite dimensional context, e.g., by the finite element method, the boundary function $\lambda_i$ can be expressed as a linear combination of suitable basis functions on $\Gamma_i$, say, $\eta^q_i(\bx), \ q=1,\ldots,N_{\Gamma_i}$, with coefficients $\Lambda^q_i$:
\begin{equation}\label{eq:lambdaRep}
 \lambda_i = \lambda_i(\bx) = \sum_{q = 1}^{N_{\Gamma_i}} \Lambda^q_i \, \eta^q_i(\bx) \, .
\end{equation}
For example, upon introducing a finite element space of continuous piecewise polynomial functions in $\Omega_i$, if $\varphi^q_i(\bx)$ are the finite element basis functions with non-null support at $\Gamma_i$, one can choose $\eta^q_i(\bx)$ to be the restriction of $\varphi^q_i(\bx)$ to $\Gamma_i$, i.e., $\eta^q_i(\bx) = \varphi^q_i(\bx)|_{\Gamma_i}$. Note that, while this is the approach used in the present work, other suitable bases can be considered on $\Gamma_i$. The dependence of the basis functions upon space is henceforth omitted, unless in the case of ambiguity.

The arbitrary coefficients $\bLambda_i = (\Lambda^1_i, \ldots, \Lambda^{N_{\Gamma_i}}_i) \in \mathcal{Q}_i$ thus become additional parameters of the local problem \eqref{eq:localProb}, with values in $ \mathcal{Q}_i = \mathcal{J}_i^1 \times \dots \times \mathcal{J}_i^{N_{\Gamma_i}}$,  where each $\mathcal{J}_i^q \subset \mathbb{R}$, $q = 1, \dots , N_{\Gamma_i}$ is a compact set.

\smallskip  

\begin{rem}
The sets of admissible boundary values introduced above need to be appropriately selected to ensure that the linear combination~\eqref{eq:lambdaRep} can approximate the trace $\lambda_i$ of the solution for all parameters $\bmu \in \mathcal{P}$. Therefore, the choice of the minimum and maximum values of $\mathcal{J}_i^q$ depends on the parameters $\bmu$.
\end{rem}

\smallskip

Although the introduction of the subdomains $\Omega_i$ in the DD procedure in Algorithm~\ref{alg:Schwarz} allows to work locally with a reduced number of spatial and parametric degrees of freedom, the parametrization of the boundary condition along the interface $\Gamma_i$ in the local problem~\eqref{eq:localProb} leads to a growth of the dimensionality of the local parametric problem, namely by introducing $N_{\Gamma_i}$ new dimensions, each associated with a coefficient $\Lambda^q_i, \ q = 1, \dots , N_{\Gamma_i}$. Unfortunately, it is well known that if $N_{\Gamma_i} \gg 1$ the solution of the local problem~\eqref{eq:localProb} with parametrized data~\eqref{eq:lambdaRep} might become unfeasible.

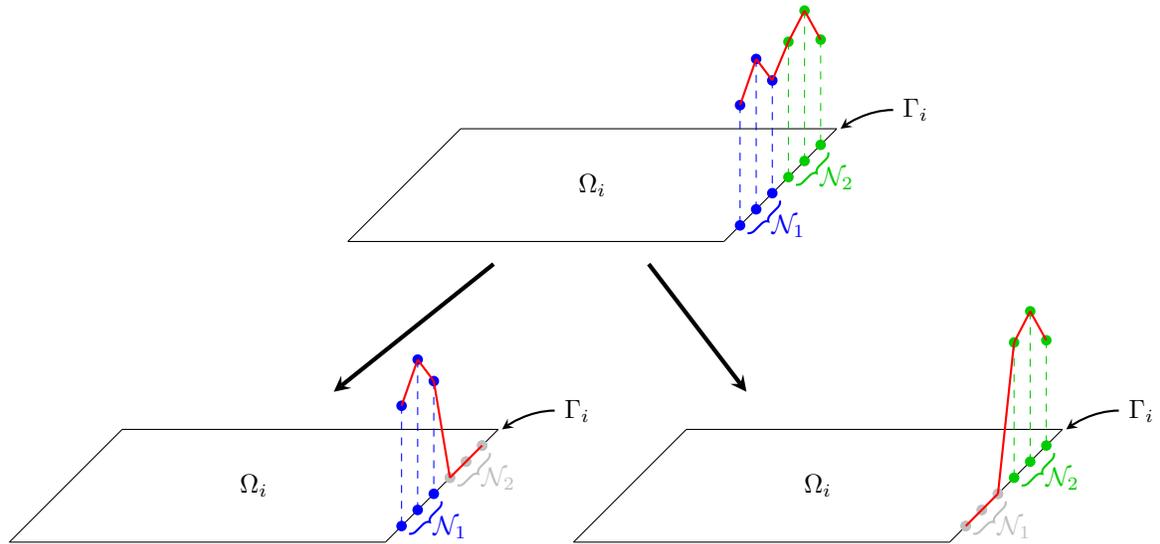
\begin{figure}[!ht]
    \centering
    \begin{tikzpicture}
        \newcommand{\intNodes}[3]{
            \draw[dashed, #3] (\longSide + #1 * \delta / 7, \ang * #1 * \delta / 7 + #2) -- (\longSide + #1 * \delta / 7, #1 * \ang * \delta / 7 );
            \fill[#3] (\longSide + #1 * \delta / 7, \ang * #1 * \delta / 7 + #2) circle (0.07 cm);
            \fill[#3] (\longSide + #1 * \delta / 7, #1 * \ang * \delta / 7 ) circle (0.07 cm);
        }

        \newcommand{\intFun}[3]{
            \draw[thick, red] (\longSide + #1 * \delta / 7, #1 * \ang * \delta / 7 + #2 ) -- (\longSide + #1 * \delta / 7 + \delta / 7, #1 * \ang * \delta / 7 +  \ang * \delta / 7 + #3);
        }

        \def\ang{1}
        \def\delta{1.5}
        \def\longSide{5}
        
        \begin{scope}
            \draw (0, 0) -- (\delta, \delta*\ang);
            \draw (0, 0) -- (\longSide, 0);
            \draw(\delta, \delta*\ang) -- (\delta + \longSide, \delta*\ang);
            \draw(\longSide, 0) -- (\delta + \longSide, \delta*\ang);
    
            \node at (\longSide / 2 + \delta / 2, \delta * \ang / 2) {$\Omega_i$};
    
            \draw[-stealth, thick] (\longSide + \delta + 0.75, \delta * \ang + 0.25) node[anchor=west] {$\Gamma_i$} parabola (\longSide + \delta + 0.05, \delta * \ang);
            
            \intNodes{1}{1.6}{blue};
            \intNodes{2}{2}{blue};
            \intNodes{3}{1.5}{blue};
    
            \intNodes{4}{1.8}{green!80!black};
            \intNodes{5}{2}{green!80!black};
            \intNodes{6}{1.4}{green!80!black};
    
            \draw [blue, decorate, decoration = {brace, mirror}, thick] (\longSide + \delta / 7 + 0.1, \ang * \delta / 7 - 0.1) -- (\longSide + 3 * \delta / 7 + 0.1, 3 * \ang * \delta / 7 - 0.1) node[pos=0.2, right=3pt, blue]{$\mathcal{N}_1$};
    
            \draw [green!80!black, decorate, decoration = {brace, mirror}, thick] (\longSide + 4 * \delta / 7 + 0.1, 4 * \ang * \delta / 7 - 0.1) -- (\longSide + 6 * \delta / 7 + 0.1, 6 * \ang * \delta / 7 - 0.1) node[pos=0.2, right=3pt, green!80!black]{$\mathcal{N}_2$};
    
            \intFun{1}{1.6}{2};
            \intFun{2}{2}{1.5};
            \intFun{3}{1.5}{1.8};
            \intFun{4}{1.8}{2};
            \intFun{5}{2}{1.4};
        \end{scope}

        \begin{scope}[shift={(-4.5, -4)}]
            \draw (0, 0) -- (\delta, \delta*\ang);
            \draw (0, 0) -- (\longSide, 0);
            \draw(\delta, \delta*\ang) -- (\delta + \longSide, \delta*\ang);
            \draw(\longSide, 0) -- (\delta + \longSide, \delta*\ang);
    
            \node at (\longSide / 2 + \delta / 2, \delta * \ang / 2) {$\Omega_i$};
    
            \draw[-stealth, thick] (\longSide + \delta + 0.75, \delta * \ang + 0.25) node[anchor=west] {$\Gamma_i$} parabola (\longSide + \delta + 0.05, \delta * \ang);
            
            \intNodes{1}{1.6}{blue};
            \intNodes{2}{2}{blue};
            \intNodes{3}{1.5}{blue};

            \intNodes{4}{0}{lightgray};
            \intNodes{5}{0}{lightgray};
            \intNodes{6}{0}{lightgray};
    
            \draw [blue, decorate, decoration = {brace, mirror}, thick] (\longSide + \delta / 7 + 0.1, \ang * \delta / 7 - 0.1) -- (\longSide + 3 * \delta / 7 + 0.1, 3 * \ang * \delta / 7 - 0.1) node[pos=0.2, right=3pt, blue]{$\mathcal{N}_1$};
    
            \draw [lightgray, decorate, decoration = {brace, mirror}, thick] (\longSide + 4 * \delta / 7 + 0.1, 4 * \ang * \delta / 7 - 0.1) -- (\longSide + 6 * \delta / 7 + 0.1, 6 * \ang * \delta / 7 - 0.1) node[pos=0.2, right=3pt, lightgray]{$\mathcal{N}_2$};
    
            \intFun{1}{1.6}{2};
            \intFun{2}{2}{1.5};
            \intFun{3}{1.5}{0};
            \intFun{4}{0}{0};
            \intFun{5}{0}{0};
        \end{scope}

        \begin{scope}[shift={(3, -4)}]
            \draw (0, 0) -- (\delta, \delta*\ang);
            \draw (0, 0) -- (\longSide, 0);
            \draw(\delta, \delta*\ang) -- (\delta + \longSide, \delta*\ang);
            \draw(\longSide, 0) -- (\delta + \longSide, \delta*\ang);
    
            \node at (\longSide / 2 + \delta / 2, \delta * \ang / 2) {$\Omega_i$};
    
            \draw[-stealth, thick] (\longSide + \delta + 0.75, \delta * \ang + 0.25) node[anchor=west] {$\Gamma_i$} parabola (\longSide + \delta + 0.05, \delta * \ang);
            
            \intNodes{1}{0}{lightgray};
            \intNodes{2}{0}{lightgray};
            \intNodes{3}{0}{lightgray};

            \intNodes{4}{1.8}{green!80!black};
            \intNodes{5}{2}{green!80!black};
            \intNodes{6}{1.4}{green!80!black};
    
            \draw [lightgray, decorate, decoration = {brace, mirror}, thick] (\longSide + \delta / 7 + 0.1, \ang * \delta / 7 - 0.1) -- (\longSide + 3 * \delta / 7 + 0.1, 3 * \ang * \delta / 7 - 0.1) node[pos=0.2, right=3pt, lightgray]{$\mathcal{N}_1$};
    
            \draw [green!80!black, decorate, decoration = {brace, mirror}, thick] (\longSide + 4 * \delta / 7 + 0.1, 4 * \ang * \delta / 7 - 0.1) -- (\longSide + 6 * \delta / 7 + 0.1, 6 * \ang * \delta / 7 - 0.1) node[pos=0.2, right=3pt, green!80!black]{$\mathcal{N}_2$};
    
            \intFun{1}{0}{0};
            \intFun{2}{0}{0};
            \intFun{3}{0}{1.8};
            \intFun{4}{1.8}{2};
            \intFun{5}{2}{1.4};
        \end{scope}
        
        \draw[- stealth, ultra thick] (1.9375, -0.3) -- (-0.1875 , -2);
        \draw[- stealth, ultra thick] (4, -0.3) -- (5.3 , -2);

    \end{tikzpicture}
    \caption{Partition of the interface nodes into two sets of active boundary parameters,  $\mathcal{N}_1$ on the left panel and $\mathcal{N}_2$ on the right panel.}
    \label{fig:activeBdryNodes}
\end{figure}

To overcome this difficulty, the linearity of the operator $L$ is exploited and the local problem~\eqref{eq:localProb} is split into a family of $N_i$ subproblems, each involving a \emph{sufficiently small} set of parameters $\mathcal{N}_i^j$, gathering the so-called \emph{active boundary parameters} (see Fig.~\ref{fig:activeBdryNodes}). More precisely,  $\{\mathcal{N}_i^j\}_{j=1,\ldots,N_i}$ denotes a disjoint partition of the set of indices $1,\ldots,N_{\Gamma_i}$ such that $\text{card}(\mathcal{N}_i^j) \ll N_{\Gamma_i}$, for all $j$.  Hence, the coefficients $\bLambda_i$ employed to characterize the trace functions can be split into subsets $\bLambda_i^j = (\Lambda_i^q)_{q \in \mathcal{N}_i^j} \in \mathcal{Q}_i^j$, with $\mathcal{Q}_i^j = \bigtimes_{q \in \mathcal{N}_i^j} \mathcal{J}_i^q \subset \mathcal{Q}_i$, and equation~\eqref{eq:lambdaRep} is rewritten as
\begin{equation}\label{eq:splittingBoundaryParameters}
\lambda_i = 
\sum_{q \in \mathcal{N}^1_i} \Lambda^q_i \, \eta^q_i + 
\sum_{q \in \mathcal{N}^2_i} \Lambda^q_i \, \eta^q_i + \ldots +
\sum_{q \in \mathcal{N}^{N_i}_i} \Lambda^q_i \, \eta^q_i .
\end{equation}

The solution of the local problem~\eqref{eq:localProb} is thus expressed in terms of both the problem parameters $\bmu \in \mathcal{P}$ and the active boundary parameters $\bLambda_i^j \in \mathcal{Q}_i^j, \ j=1,\ldots,N_i$. By linearity, for all $\bmu \in \mathcal{P}$, the resulting solution is given by
\begin{equation}\label{eq:solutionSplit}
u_i(\bmu, \bLambda_i) = u_{i,0} (\bmu) + \sum_{j=1}^{N_i} u_{i,j}(\bmu, \bLambda_i^j) ,
\end{equation}
where $u_{i,0} (\bmu)$ satisfies the equation
\begin{subequations}\label{eq:subProb}
\begin{equation}\label{eq:sourceProb}
\begin{array}{rcll}
L({u}_{i, 0}(\bmu); \bmu) &=& s_i(\bmu) & \quad \text{in } \Omega_i,\\
u_{i, 0}(\bmu) &=& g^D_i(\bmu) & \quad \text{on } \Gamma_i^D,\\
\Neum{u_{i, 0}(\bmu)}{} &=& g^N_i(\bmu) &\quad \text{on } \Gamma_i^N,\\
u_{i, 0}(\bmu) &=& 0 & \quad \text{on } \Gamma_i ,
\end{array}
\end{equation}
whereas each $u_{i, j}(\bmu, \bLambda_i^j)$, with $j=1,\ldots,N_i$, is solution of
\begin{equation}\label{eq:boundaryProbs}
\begin{array}{rcll}
L(u_{i, j}(\bmu, \bLambda_i^j); \, \bmu) &=& 0 &\quad \text{in } \Omega_i,\\
u_{i, j}(\bmu, \bLambda_i^j) &=& 0 &\quad \text{on } \Gamma_i^D,\\
\Neum{u_{i, j}(\bmu, \bLambda_i^j)}{} &=& 0 & \quad \text{on } \Gamma_i^N,\\
u_{i, j}(\bmu, \bLambda_i^j) &=& \displaystyle\sum\limits_{q \in \mathcal{N}_i^j} \Lambda^q_i\,\eta^q_i & \quad \text{on } \Gamma_i ,
\end{array}
\end{equation}
for all $\bLambda_i^j \in \mathcal{Q}_i^j$.
\end{subequations}

\subsection{Separated representation of data and local solutions}

For the sake of readability and without any loss of generality, in this section the subindex $i$ is omitted in the description of problem data in each subdomain: for instance, the Dirichlet datum $g^D(\bmu)$ is employed to seamlessly describe $g^D_i(\bmu)$ on $\Gamma_i^D$ for any subdomain $\Omega_i$.

\smallskip

In order to construct a PGD approximation of problems~\eqref{eq:subProb}, data are assumed to be given in separated form, that is,
\begin{equation}\label{eq:separatedData}
\begin{array}{c}
 \nu = \displaystyle\sum_{\ell=1}^{n_{\nu}} \xi_{\nu}^\ell(\bmu)b_{\nu}^\ell(\bx) \,, \qquad
 \balpha = \displaystyle\sum_{\ell=1}^{n_{\alpha}} \xi_{\alpha}^\ell(\bmu) \bb_{\alpha}^\ell(\bx) \,, \qquad
 \gamma = \displaystyle\sum_{\ell=1}^{n_{\gamma}} \xi_{\gamma}^\ell(\bmu) b_{\gamma}^\ell(\bx) \,, \\
 s = \displaystyle\sum_{\ell=1}^{n_s} \xi_s^\ell(\bmu) b_s^\ell(\bx) \,, \qquad
  g^N = \displaystyle\sum_{\ell=1}^{n_N}  \xi_N^\ell(\bmu) b_N^\ell(\bx) ,
\end{array}
\end{equation}
where each term of the expressions~\eqref{eq:separatedData} is the product of a function depending on the spatial coordinate $\bx$ and a function of the parameters $\bmu$. Moreover, the parametric modes are assumed to be the product of one-dimensional functions of the parameters $\mu^1,\ldots,\mu^P$, e.g.,
\begin{equation}\label{eq:vectParam}
\xi_{\nu}^\ell(\bmu) = \prod_{p=1}^P \xi_{\nu,p}^\ell(\mu^p) .
\end{equation}
Although data are not directly given in the form~\eqref{eq:separatedData}, it is possible to numerically construct a good approximation in a separated form, see~\cite{DM-MZH:15}.
In addition, the PGD rationale assumes that the solutions $u_{i, 0}$ and $u_{i, j}$ of the local subproblems can be written in separated form. 

Consider the Hilbert space
\begin{equation*}
    \mathcal{V}_i = \{w_{i} \in {H}^1(\Omega_i) \, : \, w_{i} = 0 \text{ on }\partial \Omega_i \setminus \Gamma_i^N\}\,.
\end{equation*}
The solution of the local subproblem~\eqref{eq:sourceProb} depends only on space and on the parameters $\bmu$, and it can be written as
\begin{equation}\label{eq:solProbA}
 u_{i, 0}(\bmu) = v_{i, 0}(\bmu) + G^D(\bmu) ,
\end{equation}
where $G^D(\bmu) \in H^1(\Omega_i)$ is a suitable extension of the boundary datum $g^D(\bmu)$ at $\Gamma_i^D$ such that $G^D(\bmu) = g^D(\bmu)$ at $\Gamma_i^D$ and $G^D(\bmu) = 0$ at $\Gamma_i$. Similarly to~\eqref{eq:separatedData}, also the function $G^D(\bmu)$ can be written in separated form as
\begin{equation}\label{eq:separatedDataGD}
G^D = \displaystyle\sum_{\ell=1}^{n_D} \xi_D^\ell(\bmu) b_D^\ell(\bx) \, .
\end{equation}

By construction, it follows that $v_{i, 0}(\bmu) \in \mathcal{V}_i$, for all $\bmu \in \mathcal{P}$.  Similarly, the solution of each subproblem~\eqref{eq:boundaryProbs} for $j=1,\ldots,N_i$ depends both on the parameters $\bmu$ and on the active boundary parameters $\bLambda^j_i$ at the interface, and it can be expressed as
\begin{equation}\label{eq:solProbB}
u_{i, j}(\bmu, \bLambda_i^j) = v_{i, j}(\bmu, \bLambda_i^j) + \sum_{q \in \mathcal{N}_i^j} \Lambda^q_i\,\varphi^q_i ,
\end{equation}
with $v_{i, j}(\bmu, \bLambda_i^j) \in \mathcal{V}_i$,  for all $\bmu \in \mathcal{P}$ and for all $\bLambda_i^j \in \mathcal{Q}_i^j$. 

Following the standard procedure in PGD~\cite{Chinesta:2014}, the contributions of Dirichlet boundary conditions are handled by introducing \emph{ad-hoc}, sufficiently smooth modes.  The remaining terms $v_{i, 0}(\bmu)$ and $v_{i, j}(\bmu, \bLambda_i^j)$ are computed with homogeneous Dirichlet data, under the assumption of a separated representation of all the variables, that is, $\bx$ and $\bmu$ for $v_{i,0}$, and $\bx$, $\bmu$ and $\bLambda_i^j$ for $v_{i,j}$. This yields the PGD expansions
\begin{subequations}\label{eq:nonNormalisedPGD}
\begin{align}
    v_{i,0} \approx \vpgd_{i,0} &= \sum_{m=1}^{M_0} {V}_{i,0}^m(\bx) {\phi}_{i,0}^m(\bmu) \,  ,  \label{eq:nonNormalisedPGDA} \\
    v_{i,j} \approx \vpgd_{i,j}  &= \sum_{m=1}^{M_j} {V}_{i,j}^m(\bx) {\phi}_{i,j}^m(\bmu) {\psi}_{i,j}^m(\bLambda_i^j) \, , \label{eq:nonNormalisedPGDB}
\end{align}
\end{subequations}
where ${V}_{i,0}^m$ and ${V}_{i,j}^m$ are the $m$-th spatial modes,  whereas ${\phi}_{i,0}^m$, ${\phi}_{i,j}^m$ and ${\psi}_{i,j}^m$ denote the corresponding parametric modes. It is worth noticing that the numbers of modes $M_0$ and $M_j$ are \emph{a priori} unknown and are automatically determined by a greedy procedure, see~\cite{Diez:2020:ACME}.

In the following Sects.~\ref{sect:localPGD} and \ref{sec:algebraic}, the strategy to compute \eqref{eq:nonNormalisedPGD} is presented. The result is then employed to construct $\upgd_{i, 0}(\bmu)$ and $\upgd_{i, j}(\bmu,\bLambda_i^j)$ according to~\eqref{eq:solProbA} and \eqref{eq:solProbB}. The former is the surrogate model of the data-dependent parametric problem~\eqref{eq:sourceProb}, whereas the latter are employed to define the surrogate model $\upgd_{i,\Lambda}(\bmu, \bLambda_i)$ associated with the boundary parameters, namely,
\begin{equation}\label{eq:surrogateBdry}
\upgd_{i,\Lambda}(\bmu, \bLambda_i) = \sum_{j=1}^{N_i} \upgd_{i,j}(\bmu, \bLambda_i^j) \, .
\end{equation}
Finally, the complete surrogate model for subdomain $\Omega_i$ is obtained from~\eqref{eq:solutionSplit} as
\begin{equation}\label{eq:pgdFINALsol}
\upgd_i(\bmu, \bLambda_i) = \upgd_{i,0} (\bmu) + \upgd_{i,\Lambda}(\bmu, \bLambda_i) \, .
\end{equation}

\begin{rem}
The surrogate models $\upgd_{i,j}(\bmu, \bLambda_i^j)$ feature different supports $\mathcal{Q}_i^j$ in the space of parameters $\bLambda_i^j$. 
Hence, for the summation on the right-hand side of equation~\eqref{eq:surrogateBdry} to be well-defined, each surrogate model associated with the active boundary parameters needs to be appropriately extended to have support on the entire parametric space $\mathcal{Q}_i$.
This can be straightforwardly achieved in the framework of PGD approximations by defining the modal functions for the \emph{inactive} boundary parameters $\Lambda_i^q$, $q \not\in \mathcal{N}_i^j$ to be constant and equal to $1$.
\end{rem}

\subsection{Parametric weak form of the local subproblems}
\label{sect:localPGD}
A continuous Galerkin finite element strategy is employed to construct the PGD approximations of the solutions of the local subproblems. To this end, the weak forms of the parametric problems~\eqref{eq:subProb} are first presented.

For all $v,\deV \in \mathcal{V}_i$ and for all $\bmu \in \mathcal{P}$, let $\mathcal{A}$ be the bilinear form
\begin{equation}\label{eq:bilinearA}
\mathcal{A}(v, \deV; \bmu) 
= \int_{\Omega_i} \nu(\bmu)\nabla v\cdot \nabla \deV\,d\bx 
+ \int_{\Omega_i} \balpha(\bmu) {\cdot} \nabla v \, \deV \, d\bx
+ \int_{\Omega_i} \gamma(\bmu) \, v\, \deV \,d\bx \, .
\end{equation}

The parametric weak form of problem~\eqref{eq:sourceProb} becomes: find $v_{i, 0}(\bmu) \in \mathcal{V}_i$ such that
\begin{subequations}\label{eq:pbA}
\begin{equation}\label{eq:varfA}
 \mathcal{A}(v_{i,0}(\bmu), \deV; \bmu) = \mathcal{F}_0(\deV; \bmu)\qquad \forall \deV \in \mathcal{V}_i\, \text{ and } \forall \bmu \in \mathcal{P} \, ,
\end{equation}
where
\begin{equation}\label{eq:linearA}
\mathcal{F}_0(\deV; \bmu)
= \int_{\Omega_i} s(\bmu) \deV\,d\bx
+ \int_{\Gamma_i^N} g^N(\bmu) \deV\, d\bx
- \mathcal{A}(G^D(\bmu), \deV; \bmu) 
\, .
\end{equation}
\end{subequations}

In a similar fashion,  for each problem~\eqref{eq:boundaryProbs} for $j=1,\ldots,N_i$, the parametric weak formulation is: find $v_{i, j}(\bmu, \bLambda_i^j) \in \mathcal{V}_i$ such that
\begin{subequations}\label{eq:pbB}
\begin{equation}\label{eq:varfB}
	\mathcal{A}(v_{i,j}(\bmu, \bLambda_i^j), \deV; \bmu) = \mathcal{F}_j(\deV; \bmu; \bLambda_i^j )\quad \forall \deV \in \mathcal{V}_i \, , \forall \bmu \in \mathcal{P} \, \text{ and } \forall \bLambda_i^j \in \mathcal{Q}_i^j \, ,
\end{equation}
with
\begin{equation}\label{eq:linearB}
\mathcal{F}_j(\deV; \bmu; \bLambda_i^j ) = 
- \mathcal{A} \left(\sum_{q \in \mathcal{N}_i^j} \! \Lambda^q_i\,\varphi^q_i, \deV; \bmu \right) \, .
\end{equation}
\end{subequations}

\subsection{Parametric linear systems}
\label{sec:algebraic}
Under the assumption of an affine parameter dependence of the bilinear form $\mathcal{A}$ and of the linear forms $\mathcal{F}_0$ and $\mathcal{F}_j$ (see, e.g.,~\cite{Rozza:14}), the separated approximations~\eqref{eq:nonNormalisedPGD} are constructed using a greedy approach~\cite{Chinesta:2014}. In particular,  the non-intrusive implementation provided by the encapsulated PGD solver~\cite{Diez:2020:ACME} is employed.

This approach relies on rewriting the local problems~\eqref{eq:pbA} and~\eqref{eq:pbB} in algebraic form, as parametric linear systems. To this end,  the separated representation of data~\eqref{eq:separatedData} is substituted in the parametric weak forms~\eqref{eq:varfA} and~\eqref{eq:varfB} and the unknown solutions $v_{i,0}(\bmu)$ and $v_{i,j}(\bmu, \bLambda_i^j)$ are replaced by the their corresponding PGD approximations $\vpgd_{i,0}$ and $\vpgd_{i,j}$, see~\eqref{eq:nonNormalisedPGD}.

For all $v,\deV \in \mathcal{V}_i$ and for all $\bmu \in \mathcal{P}$, let $\Apgd$ be the bilinear form
\begin{equation}\label{eq:bilinearApgd}
\begin{array}{rcl}
\Apgd(v, \deV; \bmu) 
&=& \displaystyle \sum_{\ell=1}^{n_{\nu}} \xi_{\nu}^\ell(\bmu) \int_{\Omega_i}  b_{\nu}^\ell(\bx) \nabla v\cdot \nabla \deV\,d\bx \\[3pt]
&& \displaystyle + \sum_{\ell=1}^{n_{\alpha}} \xi_{\alpha}^\ell(\bmu) \int_{\Omega_i} \bb_{\alpha}^\ell(\bx) {\cdot} \nabla v \,  \deV \, d\bx\\[3pt]
&& \displaystyle + \sum_{\ell=1}^{n_{\gamma}} \xi_{\gamma}^\ell(\bmu) \int_{\Omega_i} b_{\gamma}^\ell(\bx) v \,  \deV \,d\bx \, .
\end{array}
\end{equation}

The PGD solution $\vpgd_{i,0}$ of problem~\eqref{eq:pbA} is computed by solving the parametric equation
\begin{subequations}\label{eq:pbApgd}
\begin{equation}\label{eq:varfApgd}
 \Apgd(\vpgd_{i,0}, \deV; \bmu) = \Fpgd_0(\deV; \bmu)\qquad \forall \deV \in \mathcal{V}_i\, \text{ and } \forall \bmu \in \mathcal{P} \, ,
\end{equation}
with
\begin{equation}\label{eq:linearApgd}
\begin{array}{rcl}
\Fpgd_0(\deV; \bmu)
&=& \displaystyle \sum_{\ell=1}^{n_s} \xi_s^\ell(\bmu) \int_{\Omega_i} b_s^\ell(\bx) \deV\,d\bx \\[3pt]
&& \displaystyle + \sum_{\ell=1}^{n_N} \xi_N^\ell(\bmu) \int_{\Gamma_i^N} b_N^\ell(\bx) \deV\, d\bx \\[3pt]
&& \displaystyle - \Apgd \left(\,\sum_{\ell=1}^{n_D} \xi_D^\ell(\bmu) b_D^\ell(\bx), \deV; \bmu \right) \, .
\end{array}
\end{equation}
\end{subequations}

The PGD approximation~\eqref{eq:nonNormalisedPGDA} is constructed using a continuous Galerkin finite element discretization for each spatial mode $V_{i,0}^m(\bx)$ and a pointwise collocation approach for the parametric modes $\phi_{i,0}^m(\bmu)$. More precisely,  a finite element mesh is introduced in each subdomain $\Omega_i$ and a spatial polynomial approximation $\mathbb{Q}_r$ of degree $r \geq 1$ is selected, with basis functions $\varphi_i^n$, $n=1,\ldots,\Nfem_i$. It follows that each spatial mode can be written as
\begin{equation}\label{eq:spatialMode}
V_{i,0}^m(\bx) = \sum_{n=1}^{\Nfem_i} V_{i,0}^{m,n} \varphi_i^n(\bx) \, ,
\end{equation}
where the coefficients $V_{i,0}^{m,n}$, $n=1,\ldots,\Nfem_i$ determine the vector of spatial finite element unknowns $\mathbf{V}_{i,0}^m$. 
Therefore, the integrals appearing in the bilinear and linear forms~\eqref{eq:bilinearApgd} and~\eqref{eq:linearApgd} give rise to standard finite element matrices and vectors, appropriately weighted by means of parametric functions stemming from the separated form of data~\eqref{eq:separatedData} and~\eqref{eq:separatedDataGD}.

The resulting parametric linear system for problem~\eqref{eq:sourceProb} is
\begin{equation}\label{eq:algebraicA}
\begin{aligned}
\left(
\sum_{\ell=1}^{n_{\nu}} \xi_{\nu}^\ell(\bmu) \mat{K}_{\nu}^\ell
+ \sum_{\ell=1}^{n_{\alpha}} \xi_{\alpha}^\ell(\bmu) \mat{K}_{\alpha}^\ell
+ \sum_{\ell=1}^{n_{\gamma}} \xi_{\gamma}^\ell(\bmu) \mat{K}_{\gamma}^\ell
\right) &
\bvpgd_{i,0}(\bmu) \\
=
\sum_{\ell=1}^{n_s} \xi_s^\ell(\bmu) \mathbf{f}_s^\ell
+ \sum_{\ell=1}^{n_N} & \xi_N^\ell(\bmu) \mathbf{f}_N^\ell
+ \sum_{\ell=1}^{n_D} \xi_D^\ell(\bmu) \mathbf{f}_D^\ell
\qquad \forall \bmu \in \mathcal{P} \, ,
\end{aligned}
\end{equation}
where the PGD separated solution is defined as 
\begin{equation}\label{eq:algebraicSolA}
\bvpgd_{i,0}(\bmu) = \sum_{m=1}^{M_0} \mathbf{V}_{i,0}^m \, \phi_{i,0}^m(\bmu) \, ,
\end{equation}
whereas $ \mat{K}_{\nu}^\ell$, $ \mat{K}_{\alpha}^\ell$ and $ \mat{K}_{\gamma}^\ell$ are weighted finite element matrices stemming from the diffusion, convection and reaction term, respectively, and $\mathbf{f}_s^\ell$, $\mathbf{f}_N^\ell$ and $\mathbf{f}_D^\ell$ denote the finite element vectors accounting for the source,  Neumann and Dirichlet data, respectively.

\smallskip

The parametric linear system associated with problem~\eqref{eq:pbB} is derived with an analogous procedure.  More precisely, let $\vpgd_{i,j}$ be the solution of the parametric equation
\begin{subequations}\label{eq:pbBpgd}
\begin{equation}\label{eq:varfBpgd}
	\Apgd(\vpgd_{i,j}, \deV; \bmu) = \Fpgd_j(\deV; \bmu; \bLambda_i^j )\quad \forall \deV \in \mathcal{V}_i\, ,  \forall \bmu \in \mathcal{P} \, \text{ and } \forall \bLambda_i^j \in \mathcal{Q}_i^j \, ,
\end{equation}
with
\begin{equation}\label{eq:linearBpgd}
\Fpgd_j(\deV; \bmu;  \bLambda_i^j ) = 
- \Apgd \left( \, \sum_{q \in \mathcal{N}_i^j} \! \Lambda^q_i\,\varphi^q_i, \deV; \bmu \right) \, .
\end{equation}
\end{subequations}

The continuous Galerkin finite element discretization introduced in~\eqref{eq:spatialMode} is employed also for the spatial modes $V_{i,j}^m(\bx)$, leading to the vector of spatial unknowns $\mathbf{V}_{i,j}^m$, whereas pointwise collocation is used for the parametric modes $\phi_{i,j}^m(\bmu)$ and $\psi_{i,j}^m(\bLambda_i^j )$. Hence, the PGD approximation~\eqref{eq:nonNormalisedPGDB} is determined by computing 
\begin{equation}\label{eq:algebraicSolB}
\bvpgd_{i,j}(\bmu, \bLambda_i^j) = \sum_{m=1}^{M_j} \mathbf{V}_{i,j}^m \, \phi_{i,j}^m(\bmu) \, \psi_{i,j}^m(\bLambda_i^j )
\end{equation}
as the solution of the parametric linear system
\begin{equation}\label{eq:algebraicB}
\left(
\sum_{\ell=1}^{n_{\nu}} \xi_{\nu}^\ell(\bmu) \mat{K}_{\nu}^\ell
+ \sum_{\ell=1}^{n_{\alpha}} \xi_{\alpha}^\ell(\bmu) \mat{K}_{\alpha}^\ell
+ \sum_{\ell=1}^{n_{\gamma}} \xi_{\gamma}^\ell(\bmu) \mat{K}_{\gamma}^\ell
\right)
\bvpgd_{i,j}(\bmu,\bLambda_i^j) 
=
\sum_{q \in \mathcal{N}_i^j} \Lambda^q_i \mathbf{f}_{\Lambda}^q \, ,
\end{equation}
for any value of $\bmu \in \mathcal{P}$ and $\bLambda_i^j \in \mathcal{Q}_i^j$.  In equation~\eqref{eq:algebraicB}, the vector $\mathbf{f}_{\Lambda}^q$ stems from imposing the parametric Dirichlet boundary condition at the interface $\Gamma_i$ in equation~\eqref{eq:boundaryProbs} within the finite element setting.

\smallskip

The encapsulated PGD library~\cite{Diez:2020:ACME} is utilized to solve equations~\eqref{eq:algebraicA} and \eqref{eq:algebraicB}. Technical details on the setup of problems~\eqref{eq:algebraicA} and~\eqref{eq:algebraicB} in the encapsulated PGD framework for a sample test case are presented in Appendix~\ref{append:encapsulatedPGD}.

\smallskip

\begin{rem}
For the case of two subdomains, the cost of the offline phase stems from the computation of two surrogate models accounting for the data-dependent parametric problems~\eqref{eq:sourceProb} and $N_1+N_2$ surrogate models related to problems~\eqref{eq:boundaryProbs} with active boundary parameters.
It is worth noticing that all PGD approximations mentioned above are independent from one another and can be efficiently computed in parallel. 
Moreover, the computational effort during the offline phase can be further reduced, e.g.,  by identifying a reference subdomain where local surrogate models are computed before being suitably mapped to the physical subdomains of the problems under consideration, as demonstrated in the example in Sect.~\ref{sec:testPatera}.
\end{rem}

\section{Surrogate-based overlapping Schwarz method}
\label{sec:online}

In this section,  an efficient strategy to construct the global solution of problem~\eqref{eq:globalProb} for a fixed set of parametric values $\bar{\bmu} \in \mathcal{P}$ is presented.
The goal is to devise a procedure, suitable for real time execution, to appropriately glue the parametric solutions of the local subproblems, thus drastically reducing the cost of the overall DD algorithm.
To this end, the overlapping Schwarz method presented in Sect.~\ref{sec:schwarzAlg} is adapted to exploit the local PGD surrogate models constructed in Sect.~\ref{sec:offline}.

\begin{rem}
To easily perform the coupling between subdomains $\Omega_1$ and $\Omega_2$ in the online phase and to avoid expensive interpolation procedures among different grids, the meshes used for the spatial discretization of the local subdomains are assumed to be conforming with the interfaces (i.e., the interfaces do not cut through any elements of the meshes) and to coincide in the overlapping region $\Omega_{12}$.
\end{rem}

First, note that by construction (see~\eqref{eq:lambdaRep}),  the vector $\mathbf{u}_{\Gamma_i}\!(\bar{\bmu})$ of the nodal values of the solution at $\Gamma_i$ corresponds to the vector of parameters $\bLambda_i$. Hence, for any $\bar{\bmu} \in \mathcal{P}$, equation~\eqref{eq:SchwarzAlg} can be rewritten as
\begin{equation}\label{eq:SchwarzAlgLambda}
\mat{I}_{\Gamma_j}\bLambda_j
=
\mat{R}_{\Omega_i\to\Gamma_j} \mat{A}_{\Omega_i}^{-1} \mathbf{f}_{\Omega_i}\!(\bar{\bmu}) 
+ \mat{R}_{\Omega_i\to\Gamma_j} \mat{A}_{\Omega_i}^{-1} (-\mat{A}_{\Gamma_i}\bLambda_i )\, ,
\end{equation}
where $\mat{A}_{\Omega_i}^{-1} \mathbf{f}_{\Omega_i}\!(\bar{\bmu})$ corresponds to the solution of problem~\eqref{eq:sourceProb} and $\mat{A}_{\Omega_i}^{-1} (-\mat{A}_{\Gamma_i}\bLambda_i)$ denotes the solution of problem~\eqref{eq:boundaryProbs}.
Exploiting the local surrogate models constructed in the offline phase, equation~\eqref{eq:SchwarzAlgLambda} reduces to
\begin{equation}\label{eq:SchwarzAlgPGD}
\mat{I}_{\Gamma_j}\bLambda_j
=
\mat{R}_{\Omega_i\to\Gamma_j} \bupgd_{i,0}\!(\bar{\bmu}) 
+ \mat{R}_{\Omega_i\to\Gamma_j} \bupgd_{i,\Lambda}\!(\bar{\bmu},\bLambda_i)  \, ,
\end{equation}
where $\bupgd_{i,0}\!(\bar{\bmu})$ and $\bupgd_{i,\Lambda}\!(\bar{\bmu},\bLambda_i)$ denote the vectors of the nodal values of the PGD solutions $\upgd_{i,0}$ and $\upgd_{i,\Lambda}$, respectively, evaluated for the target values $\bar{\bmu} \in \mathcal{P}$ and $\bLambda_i \in \mathcal{Q}_i$.

Let $\mat{A}_{i,\Lambda_i}^\texttt{PGD}$ be the local PGD operator 
\begin{equation}\label{eq:PGDlambdaOper}
\mat{A}_{i,\Lambda_i}^\texttt{PGD} : \bLambda_i \to \bupgd_{i,\Lambda}(\bar{\bmu},\bLambda_i) \, 
\end{equation}
that, given a set of boundary parameters $ \bLambda_i$, returns the nodal values of the PGD surrogate model $\upgd_{i,\Lambda}$ of problem~\eqref{eq:boundaryProbs} for the set of parameters $\bar{\bmu}$. Hence, equation \eqref{eq:SchwarzAlgPGD} can be rewritten as
\begin{equation}\label{eq:SchwarzAlgPGD_Op}
\mat{I}_{\Gamma_j}\bLambda_j
=
\mat{R}_{\Omega_i\to\Gamma_j} \bupgd_{i,0}\!(\bar{\bmu}) 
+ \mat{R}_{\Omega_i\to\Gamma_j} \mat{A}_{i, \Lambda_i}^\texttt{PGD} \bLambda_i  \, ,    
\end{equation}
and the surrogate-based overlapping Schwarz method is finally obtained by rewriting the interface system~\eqref{eq:systemSchwarzInterface} as
\begin{equation}\label{eq:systemSchwarzInterfacePGD}
\begin{pmatrix}
\mat{I}_{\Gamma_1} & -\mat{R}_{\Omega_2\to\Gamma_1} \mat{A}_{2,\Lambda_2}^\texttt{PGD} \\[3pt]
-\mat{R}_{\Omega_1\to\Gamma_2} \mat{A}_{1,\Lambda_1}^\texttt{PGD} & \mat{I}_{\Gamma_2}
\end{pmatrix}
\begin{pmatrix}
\bLambda_1 \\[3pt]
\bLambda_2
\end{pmatrix}
=
\begin{pmatrix}
\mat{R}_{\Omega_2\to\Gamma_1} \bupgd_{2,0}\!(\bar{\bmu}) \\[3pt]
\mat{R}_{\Omega_1\to\Gamma_2} \bupgd_{1,0}\!(\bar{\bmu})
\end{pmatrix} \, .
\end{equation}

Therefore, the online phase of the method consists of an iterative strategy to solve equation~\eqref{eq:systemSchwarzInterfacePGD}, e.g., by GMRES.  At convergence, say, at iteration $k=k^*$,  the approximation of the solution of the global problem~\eqref{eq:globalProb} for $\bar{\bmu} \in \mathcal{P}$ is thus given by
\begin{equation}\label{eq:globalSolnA}
	\bupgd(\bar{\bmu}) = 
	\begin{cases}
		\bupgd_{1,0}(\bar{\bmu}) + 
		\bupgd_{1,\Lambda}(\bar{\bmu},\bLambda_1^{(k^*)} ) \quad \text{in } \Omega_1,\\
		\noalign{\vskip5pt}
		\bupgd_{2,0}(\bar{\bmu}) + 
		\bupgd_{2,\Lambda}(\bar{\bmu},\bLambda_2^{(k^*)} ) \quad \text{in } \Omega_2 \setminus \Omega_{12}.
	\end{cases}
\end{equation}

It is worth noticing that the values $\bLambda_i^{(k)}$ computed by GMRES iterations to solve problem~\eqref{eq:systemSchwarzInterfacePGD} may not coincide with any of the values obtained from the discretization of the parametric domain $\mathcal{Q}_i$.  
If this is the case, the solution $\bupgd_{i,\Lambda}(\bar{\bmu},\bLambda_i^{(k)})$ provided by the operator $\mat{A}_{i,\Lambda}^\texttt{PGD}$ is obtained by performing a linear interpolation of the parametric modes depending on $\bLambda_i$ and associated with the available values closest to $\bLambda_i^{(k)}$.

\begin{rem}
Following Algorithm~\ref{alg:Schwarz}, at the beginning of the online phase,  an instance $\bar{\bmu}$ is selected in the set of parameters $\mathcal{P}$.
This is not strictly necessary and the described algorithm can be adapted to handle arbitrary parameters $\bmu$. 
In the latter case, at the end of the online phase, one would obtain a global surrogate model that represents a family of solutions of the global problem~\eqref{eq:globalProb} depending on $\bmu \in \mathcal{P}$, instead of an instance of such model for $\bmu = \bar{\bmu}$. 
\end{rem}

\section{Numerical tests}
\label{sec:results}

In this section, some numerical tests\footnote{The numerical results presented in this section have been obtained using a PC with CPU Intel$^{\mbox{\tiny{\textregistered}}}$ Core\texttrademark\; i5-11400 @ 2.60GHz and 8GB RAM.} are presented to assess the performance of the proposed DD-PGD method, considering three elliptic problems: a diffusion problem with synthetic solution in two subdomains (Sect.~\ref{sec:testAnalytical}), a convection-diffusion problem with parametric geometry and two subdomains in Sect.~\ref{sec:testRozza}, and a \rev{thermal problem featuring multiple subdomains, each with a different bulk conductivity, in Sect.~\ref{sec:testPatera}.}

The local parametric subproblems are solved using the encapsulated PGD toolbox~\cite{Diez:2020:ACME}. Unless otherwise specified, a tolerance $\varepsilon = 10^{-4}$ is selected for the PGD enrichment process and the redundant information is then eliminated by the PGD compression algorithm~\cite{DM-MZH:15} with tolerance $\varepsilon^\star = 10^{-3}$. 

In the offline phase, the meshes employed for the spatial discretization are problem-dependent and they are specified in the corresponding sections. It is worth noticing that, for all considered numerical tests, the local meshes in the spatial subdomains are conforming with the interfaces (i.e., the interfaces do not cut through any elements of the meshes) and they coincide in the overlapping regions. The one-dimensional parametric intervals $\mathcal{I}^p$, $p=1,\ldots,P$ and $\mathcal{J}_i^q$, $q=1,\ldots,N_{\Gamma}$ are discretized using uniform elements.
In the online phase, the interface system is solved using GMRES~\cite{Saad:1986:SISSC}, with tolerance $10^{-6}$ on the relative residual. 

\rev{
\begin{rem}
The number of active boundary parameters to be employed in each subproblem is case-dependent.
Whilst a large set of active boundary parameters reduces the number of subproblems to be solved, the performance of the PGD-ROM algorithm can be negatively affected in the presence of high-dimensional problems. Indeed, the implementation used to generate the results in this section requires fewer modes to converge when fewer parameters are considered at the interface, thus leading to computational gains when local problems are solved in parallel.
In this context, subproblems with only one active boundary parameter can also be devised. This case is equivalent to computing local surrogate models setting unitary Dirichlet boundary conditions at each node of the interface, and appropriately scaling and combining the resulting PGD solutions during the online phase. The two approaches provide comparable solutions in terms of the accuracy of the global surrogate model, the number of PGD modes and the overall computing time, and they are not reported here for brevity.
In order to reduce the number of local surrogate models to be determined whilst also avoiding high-dimensional problems, in the rest of this section, problems of type~\eqref{eq:boundaryProbs} are solved with at most three active boundary parameters. 
However, further research is required to understand how to optimally choose the sets of active boundary parameters in the DD-PGD framework (see, e.g., optimal port spaces in RBE~\cite{Patera-SP-16}), in order to balance the accuracy of the local surrogate models, the computational complexity of each subproblem and their number. 
\end{rem}
}

\subsection{Diffusion problem with synthetic solution in two subdomains}
\label{sec:testAnalytical}

This model problem aims to assess the convergence and accuracy of the proposed method. Let $\Omega = (0, 2) \times (0, 1)$ be the spatial domain and $\mu \in \mathcal{P} = [1,50]$ be a scalar parameter defining a space-dependent conductivity coefficient. The parametric Poisson equation under analysis is: for all $\mu \in \mathcal{P}$, find $u(\mu)$ such that 
\begin{equation}
    \label{eq:numTestsModel}
	\begin{array}{rcll}
		- \nabla \cdot ((1 + \mu x) \nabla u(\mu)) &=& s(\mu) & \quad \text{in } \Omega \, ,\\
		u(\mu) &=& 0 & \quad \text{on } \partial\Omega \, ,\\
	\end{array}
\end{equation}
where the source term $s(\mu)$, detailed in Appendix~\ref{append:Setup}, is selected such that the analytical solution of the problem is ${u}(\mu) = {u_{\text{ex}}}(\mu)$ with
\begin{equation}
    \label{eq:numTests1analytical}
    {u_{\text{ex}}}(\mu) := \sin(2\pi x)\sin(2\pi y) + \frac{\mu}{2}xy(y - 1)(x - 2)\,.
\end{equation}

For the spatial discretization, continuous $\mathbb{Q}_1$ Lagrangian finite elements are employed on a structured grid with local mesh size $\hX = 5 \times 10^{-2}$. The domain $\Omega$ is split into the two overlapping subdomains $\Omega_1 = (0, 1 + n\hX) \times (0, 1)$ and $\Omega_2 = (1 - n\hX, 2) \times (0, 1)$, where $n \in \mathbb{N}$ is a positive value that controls the size $2n\hX$ of the overlap. For instance, if $n = 1$ the width is equal to $2\hX$, that is, the overlap region contains two layers of elements. The mesh of each subdomain contains a total of 420 quadrilateral elements.

Following the definition of the spatial mesh, 19 nodes are considered as auxiliary boundary parameters $\Lambda_i^q$ at the interfaces $\Gamma_1 = \{1+n\hX\}\times(0,1)$ and $\Gamma_2 = \{1-n\hX\}\times(0,1)$,  excluding the nodes at $\Gamma_1 \cap \partial\Omega$ and $\Gamma_2 \cap \partial\Omega$. Since $|u(\mu)|<10$ for all $\mu\in\mathcal{P}$ and for all $(x,y)\in\Omega$, the range of admissible values for $\Lambda_i^q$ is set to $\mathcal{J}_i^q = [-10, 10]$, for $i=1,2$ and for $q=1,\ldots,19$. All parametric domains (i.e., $\mathcal{P}$ and $\mathcal{J}_i^q$, for $i=1,2$ and for $q=1,\ldots,19$) are discretized by equally-spaced nodes, with grid spacing $\hMu = \hLam = 10^{-3}$, leading to $\Nmu = 4.9 \times 10^{4}$ nodes in $\mathcal{P}$ and $\Nlam = 2 \times 10^{4}$ nodes in $\mathcal{J}_i^q$, for all $q$.

Each local problem is thus 22-dimensional, featuring two spatial dimensions, one scalar physical parameter $\mu$ and 19 boundary parameters $\Lambda_i^q$ at the interface. To reduce the dimensionality of the problem, the interface boundary parameters are partitioned into sets of at most \rev{3} active parameters. Hence, for each subdomain, one problem of type~\eqref{eq:sourceProb} of dimension 3 and 6 problems of type~\eqref{eq:boundaryProbs}, with dimension at most \rev{6}, are identified.

\begin{figure}[h!]
    \centering
    \includegraphics{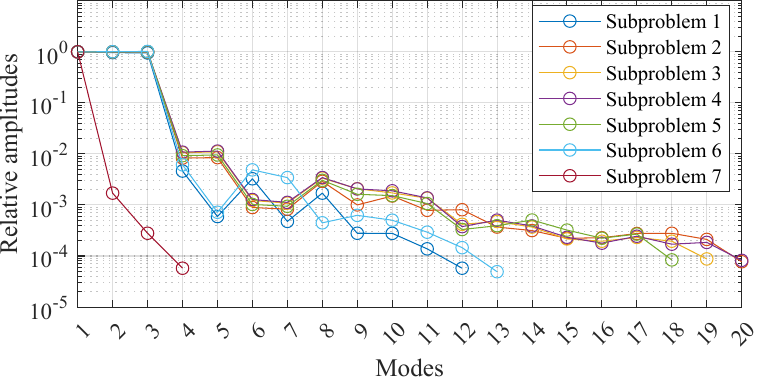}
    \includegraphics{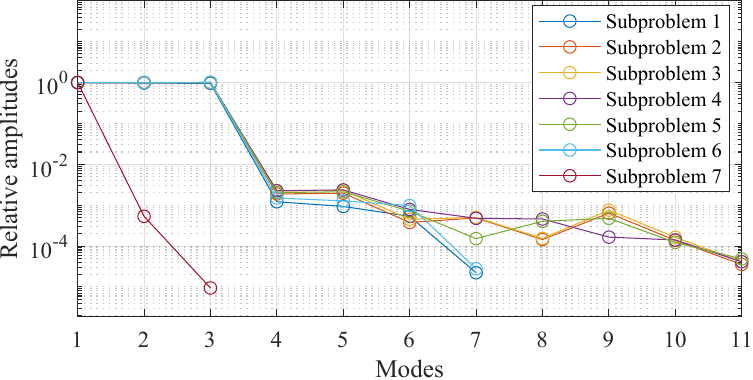}
    \caption{\rev{Relative amplitudes of the PGD modes, scaled by the amplitude of the first mode, for subdomain 1 (top) and subdomain 2 (bottom).}}
    \label{fig:PGDAmps}
\end{figure}

The PGD solutions of the local parametric subproblems are formed by at most $116$ modes, and the total computational time of the offline phase is approximately \SI{98}{s}.
\rev{More precisely, the convergence of the PGD enrichment is described in Fig.~\ref{fig:PGDAmps}, reporting the amplitude of the computed modes scaled by the amplitude of the first mode, for all subproblems in each subdomain. It is worth noticing that, in each subdomain, subproblem 7  only contains 2 active boundary parameters, whence the significant reduction of required modes with respect to the remaining subproblems. Moreover, the}
PGD compression reduces the maximum number of modes required for a local problem to $59$, with an extra execution time in the offline phase of \SI{13}{s}.
\rev{
For the sake of brevity, the convergence of the PGD enrichment is only presented for this test case, while for the following problems only the total number of modes required in each subdomain is detailed and the presented results only report the cost of the compressed local PGD solutions.
}

\begin{figure}[h!]
    \centering
    \includegraphics{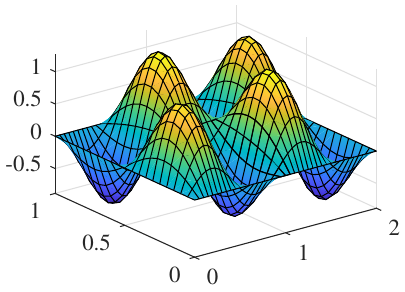}
    \includegraphics{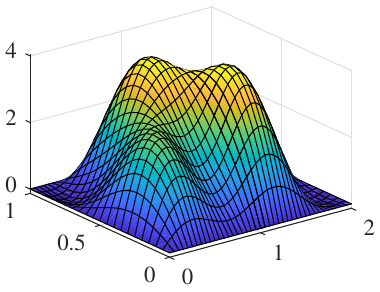}
    \caption{Global DD-PGD solution for $\mu = 3$ (left) and $\mu = 30$ (right) and overlap of width $2\hX=10^{-1}$.}
    \label{fig:testAnalyticalSolutions}
\end{figure}

Setting the maximum dimension of the Krylov subspace to 6 in the GMRES algorithm, the DD-PGD strategy is considered to evaluate the global solution for $\mu=3$ and $\mu=30$. The PGD solutions, computed with overlap of width $2\hX$ (i.e., $n=1$), are plotted in Fig.~\ref{fig:testAnalyticalSolutions}.

\begin{figure}[h!]
    \centering
        \includegraphics{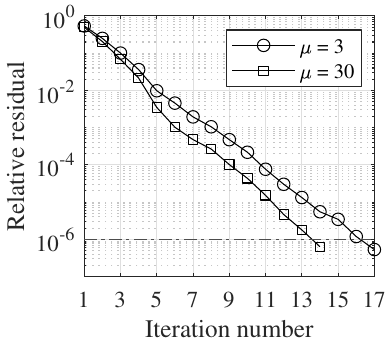}
        \hspace{1cm}
        \includegraphics{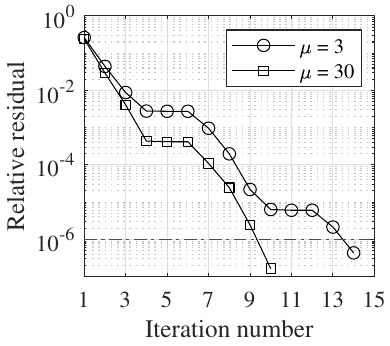}
        \caption{GMRES relative residual as a function of the number of iterations for $\mu=3,\,30$ and overlap of width $2\hX=10^{-1}$ (left) and $6\hX=3 \times 10^{-1}$ (right).}
        \label{fig:testAnalyticalResidual}
\end{figure}

The convergence of the GMRES relative residuals as a function of the iterations is shown in Fig.~\ref{fig:testAnalyticalResidual}, for different widths of the overlapping region. 
The results display that fewer iterations are needed when the width of the overlap increases from $2\hX$ (i.e., $n=1$) to $6\hX$ (i.e., $n=3$), a typical behaviour for the one-level Schwarz method, see, e.g.,~\cite{Smith:1996}.

For the overlap width $2\hX$ and for $\mu = 3$, the online phase of the DD-PGD approach requires 17 iterations to achieve the user-defined tolerance (see Fig.~\ref{fig:testAnalyticalResidual}, left), with a computational time of approximately \SI[parse-numbers=false]{4 \times 10^{-2}}{s}. To assess the cost of this procedure, a standard DD-FEM approach is executed using the same overlap width, for 50 random values of the parameter $\mu$. The high-fidelity DD-FEM requires an average of 9 iterations to converge to the specified tolerance, with a mean computing time of \SI{1.85}{s}. 
Although the DD-PGD requires a larger number of GMRES iterations to achieve the prescribed tolerance, likely as a consequence of the loss of information occurring during the ROM construction, it is worth noticing that the method is still approximately 46 times faster than the standard DD-FEM procedure, owing to the real-time evaluation of the local surrogate solutions.

Next, the accuracy of the proposed methodology is assessed for the previously selected values of the parameter $\mu$. More precisely, Table~\ref{tab:numTest1Errors} reports the relative error, measured in the $L^2(\Omega)$ norm, between the analytical solution $u_{\text{ex}}$ and the DD-PGD solution $\upgd$, the monolithic PGD solution $\upgd_{\Omega}$ computed on the entire domain $\Omega$ without partitioning it into $\Omega_1$ and $\Omega_2$ and the high-fidelity monolithic finite element solution $u^h_{\Omega}$. 
The results clearly display that the DD-PGD strategy provides a solution with global accuracy comparable both to the monolithic ROM and the high-fidelity solution. In particular, it is worth remarking that the monolithic PGD surrogate model $\upgd_{\Omega}$ involves the solution of spatial problems almost twice as large as the ones appearing in the DD-PGD framework, with an overall computing time of \SI[parse-numbers=false]{4.5 \times 10^{-1}}{s}, contrary to the \SI[parse-numbers=false]{4 \times 10^{-2}}{s} of the DD-PGD approach.

\begin{figure}[h!]
    \centering
\includegraphics{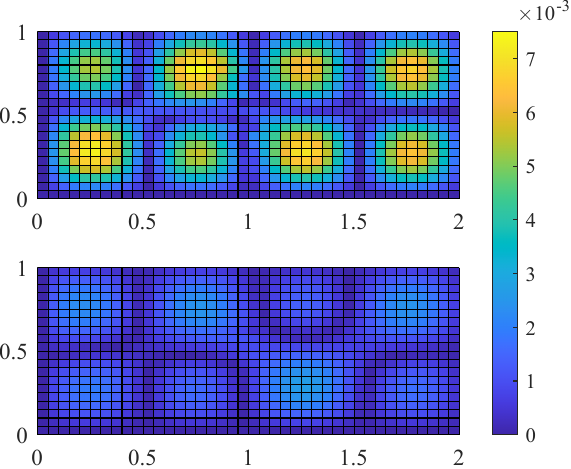}
    \caption{Map of the scaled nodal error $|\upgd(\mu) - u_{\text{ex}}(\mu)| / \max_{\Omega} |u_{\text{ex}}(\mu) |$ in the domain $\Omega$, for $\mu=3$ (top) and $\mu=30$ (bottom).}
    \label{fig:testAnalyticalError}
\end{figure}

Finally, Fig.~\ref{fig:testAnalyticalError} reports the map of the scaled nodal error $|\upgd(\mu) - u_{\text{ex}}(\mu)| / \max_{\Omega} |u_{\text{ex}}(\mu)|$ in the entire domain $\Omega$ for different values of $\mu$, indicating that no loss of accuracy is experienced by the DD-PGD solution in the neighbourhood of the overlapping region.

\begin{table}[h!]
    \centering
    \begin{tabular}{cccc}
        $\mu$ & $\frac{\|\upgd(\mu) - u_{\text{ex}}(\mu)\|_{L^2(\Omega)}}{\|u_{\text{ex}}(\mu)\|_{L^2(\Omega)}}$ &  $\frac{\|\upgd_{\Omega}(\mu) - u_{\text{ex}}(\mu)\|_{L^2(\Omega)}}{\|u_{\text{ex}}(\mu)\|_{L^2(\Omega)}}$ & $\frac{\|u^h_{\Omega}(\mu) - u_{\text{ex}}(\mu)\|_{L^2(\Omega)}}{\|u_{\text{ex}}(\mu)\|_{L^2(\Omega)}}$ \\
        \hline
         3 & $9.27 \times 10^{-3}$ & $9.00 \times 10^{-3}$ &  $9.07 \times 10^{-3}$ \\
         30 & $3.26 \times 10^{-3}$ & $3.27 \times 10^{-3}$  &  $3.27 \times 10^{-3}$ \\
         \hline
    \end{tabular}
    \caption{Relative error in $L^2(\Omega)$ norm between the analytical solution $u_{\text{ex}}(\mu)$ and the DD-PGD solution $\upgd(\mu)$, the monolithic PGD solution $\upgd_{\Omega}(\mu)$ and the monolithic finite element solution $u^h_{\Omega}(\mu)$.}
    \label{tab:numTest1Errors}
\end{table}

\subsection{Poiseuille–Graetz flow in a geometrically parametrized domain}
\label{sec:testRozza}

In this section, a convection-diffusion equation in a parametrized domain, presented in~\cite{Pacciarini:2014:CMAME}, is studied. 
This benchmark describes a channel, with walls maintained at different temperatures, in which heat convection and conduction phenomena are combined. 

\begin{figure}[h!]
    \centering
    \begin{tikzpicture}
        \draw (0, 0) node[anchor = north]{$(0,0)$} -- (1.5, 0) node[anchor = north]{$\Gamma^{D, 1}$} -- (3, 0) node[anchor = north east, shift={(0.25, 0)}]{$(1, 0)$} ;
        \draw (0, 0) -- (0, 1.5) node[anchor = east]{$\Gamma^{D, 1}$} -- (0, 3)  node[anchor = south]{$(0, 1)$};
        \draw (0, 3) -- (1.5, 3) node[anchor = south]{$\Gamma^{D, 1}$} -- (3, 3);
        \draw[<->, thick, dashed] (0, 1.2) -- (3.2, 1.2);

        \node at (1.5, 1.5) {$\Omega_1$};

        \draw[dotted] (3, 0) -- (3, 3);
        \draw[dotted] (3.2, 0) -- (3.2, 3) node[anchor = south west, shift={(-0.2, 0)}]{$(1.05, 1)$};

        \draw[ultra thick] (3, 0) -- (5.25, 0) node[anchor = north]{$\Gamma^{D, 2}$} --(7.5, 0) node[anchor = north]{$(1 + \mu_2, 0)$};
        \draw[ultra thick] (3, 3) -- (5.25, 3) node[anchor = south]{$\Gamma^{D, 2}$} --(7.5, 3) node[anchor = south]{$(1 + \mu_2, 1)$};
        \draw[dashed] (7.5, 0) -- (7.5, 1.5) node[anchor = west]{$\Gamma^N$} -- (7.5, 3);
        \draw[<->, thick, dashed] (3, 1.8) -- (7.5, 1.8);

        \node at (5.25, 1.5) {$\Omega_2$};
    \end{tikzpicture}
    \caption{Partition of the parametrized domain for the Poiseuille-Graetz flow problem.}
    \label{fig:PGflowDom}
\end{figure}
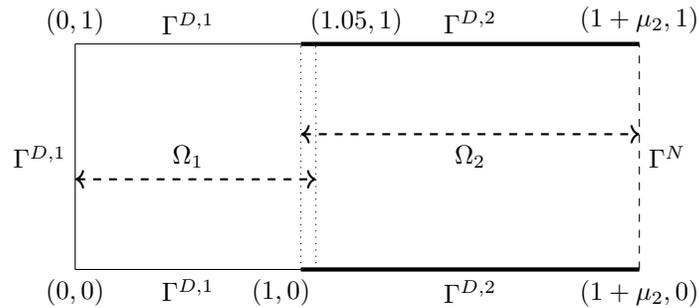

The test features two parameters: a physical one, $\mu_1$, describing the inverse of the diffusion coefficient, and a geometric one, $\mu_2$, controlling the size of the domain $\Omega = (0, 1 + \mu_2) \times (0, 1)$ (see Fig.~\ref{fig:PGflowDom}). 

For all $\bmu = (\mu_1, \mu_2) \in \mathcal{P}$, with $\mathcal{P} = [10^4, 2 \times 10^4] \times [0.5, 4]$, a solution $u(\bmu)$ is sought to fulfil the problem
\begin{equation}\label{eq:pgProblem}
    \begin{array}{rcll}
    \displaystyle
         -\frac{1}{\mu_1} \,\Delta u(\bmu) + \balpha \cdot \nabla u(\bmu) &=& 0 &\quad \text{in } \Omega(\bmu)\, ,\\
         u(\bmu) &=& 0 &\quad \text{on } \Gamma^{D, 1}(\bmu)\, ,\\
         u(\bmu) &=& 1 &\quad \text{on } \Gamma^{D, 2}(\bmu)\, ,\\
         \displaystyle\frac{1}{\mu_1}\nabla u(\bmu)\cdot \bn(\bmu) &=& 0 &\quad \text{on } \Gamma^N(\bmu) \, ,
    \end{array}
\end{equation}
where $\balpha = (\alpha_1,\alpha_2)^T = (4y(1 - y) , 0 )^T$ is a parameter-independent, horizontal velocity field and the boundaries of the parametric domain are such that
\begin{align*}
    \Gamma^{D, 1}(\bmu) &= \{0\} \times [0, 1] \cup [0, 1] \times \{0\} \cup [0, 1] \times \{1\} \,,\\
    \Gamma^{D, 2}(\bmu) &= [1, 1 + \mu_2] \times \{0\} \cup [1, 1 + \mu_2] \times \{1\}\, ,\\
    \Gamma^N(\bmu) &= \{1 + \mu_2\} \times [0, 1] \, .
\end{align*}

The domain $\Omega(\bmu)$ is partitioned into two overlapping regions, the parameter-independent subdomain $\Omega_1 = [0, 1.05] \times [0, 1]$ and the parametric subdomain $\Omega_2(\bmu) = [1, 1 + \mu_2] \times [0, 1]$. 
In order to devise a surrogate model for the geometrically parametrized subdomain $\Omega_2(\bmu)$, a formulation based on a reference domain configuration~\cite{Ammar-AHCCL-14} is employed. This is achieved by introducing a parameter-independent reference subdomain $\hOmega_2 = [0,1] \times [0,1]$ and an appropriate parametric mapping $\Map : \hOmega_2 \times \mathcal{I}^2 \rightarrow \Omega_2(\bmu)$ to transform the fixed reference subdomain into the parametric physical one.

\begin{figure}[h!]
    \centering
    \includegraphics[]{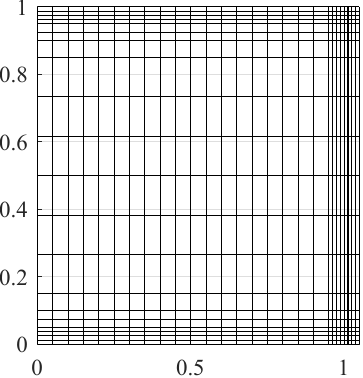}
    \hspace{1cm}
    \includegraphics[]{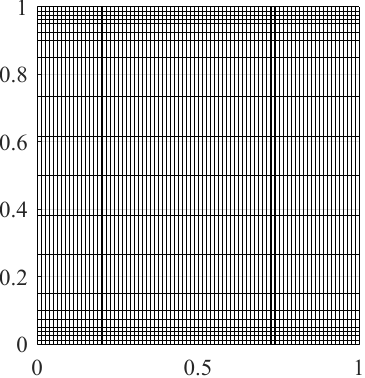}
    \caption{Computational meshes for subdomains $\Omega_1$ (left) and $\hOmega_2$ (right) in the Poiseuille-Graetz flow problem.}
    \label{fig:meshSize}
\end{figure}

Structured grids containing 540 and 1,600 quadrilaterals are defined in $\Omega_1$ and $\hOmega_2$, respectively, with non-uniform mesh size near the top and bottom walls, as detailed in Fig.~\ref{fig:meshSize}. 

\begin{figure}[h!]
    \centering
    \includegraphics[]{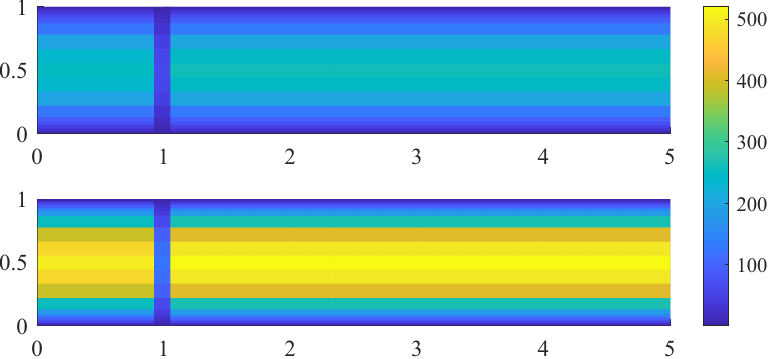}
    \caption{Local P\'eclet number in the horizontal direction computed element-by-element for the case of $\mu_1 = 10^4$ (top) and $\mu_1 = 2 \times 10^4$ (bottom).}
    \label{fig:Peclet}
\end{figure}

Since the second component of the convective velocity is null, the local P\'eclet number in the horizontal direction (i.e., $\mathrm{Pe} = |\alpha_1| h_{x,1} \mu_1 /2$, $h_{x,1}$ being the local mesh size in the horizontal direction) is displayed in Fig.~\ref{fig:Peclet} for different values of the diffusion coefficient $1/\mu_1$. For all the considered values of the parameter $\mu_1$, problem~\eqref{eq:pgProblem} is convection-dominated. Hence, following~\cite{Huerta-GCCDH-13}, a continuous $\mathbb{Q}_1$ Lagrangian finite element formulation with streamline upwind Petrov-Galerkin (SUPG) stabilization is implemented in the PGD surrogate model. Details of the SUPG formulation on a parameter-independent reference domain are provided in Appendix~\ref{append:Graetz}.

The interface boundary parameter intervals are selected as $\mathcal{J}_i^q = [-5, 5]$, for $i=1,2$ and for all $q$. All parametric intervals are discretized using equidistributed nodes, with spacing $h_{\mu_1} = 10^{-1}$ and $h_{\mu_2} = h_{\Lambda} = 10^{-3}$, leading to $N^{\mu_1} = 10^{5}$, $N^{\mu_2} = 3.5 \times 10^{3}$ and $\Nlam = 10^{4}$ parametric unknowns in $\mathcal{I}^1$, $\mathcal{I}^2$ and $\mathcal{J}_i^q$, respectively.

For each subdomain, one subproblem of type~\eqref{eq:sourceProb} is solved. Moreover, 7 subproblems of type~\eqref{eq:boundaryProbs} are solved in $\Omega_1$, with three active interface boundary parameters, whereas in subdomain $\hOmega_2$, 10 such subproblems are formulated, with at most two active interface boundary parameters. This choice follows from problems in subdomain $\hOmega_2$ being higher dimensional, due to the presence of the geometric parameter $\mu_2$.

The DD-PGD strategy constructs local surrogate models featuring 56 modes in $\Omega_1$ and 619 modes in subdomain $\hOmega_2$. After compression, 38 and 164 modes are obtained for $\Omega_1$ and $\hOmega_2$, respectively, for a total CPU time of approximately \SI[parse-numbers=false]{6.88 \times 10^{3}}{s}.
It is worth noticing that the imbalance in the computational effort in $\Omega_1$ and $\hOmega_2$ stems from the second subdomain featuring two parameters, with one of them (i.e., $\mu_2$) controlling the geometric transformation of the domain.
Indeed, it is well known that parametric variations of the geometry and the concurrent presence of multiple parameters cause a significant increase in the complexity of parametric problems, thus in the number of terms required in the PGD approximations~\cite{Giacomini-GBSH-21}.
In this context, the DD-PGD strategy is particularly appealing as it allows to separately treat the portion of the domain with multiple parameters and complex dynamics by computing a large number of local modes (i.e., 619 in $\hOmega_2$), while only 56 terms need to be determined in the region $\Omega_1$.

\begin{figure}[h!]
    \centering
    \includegraphics[]{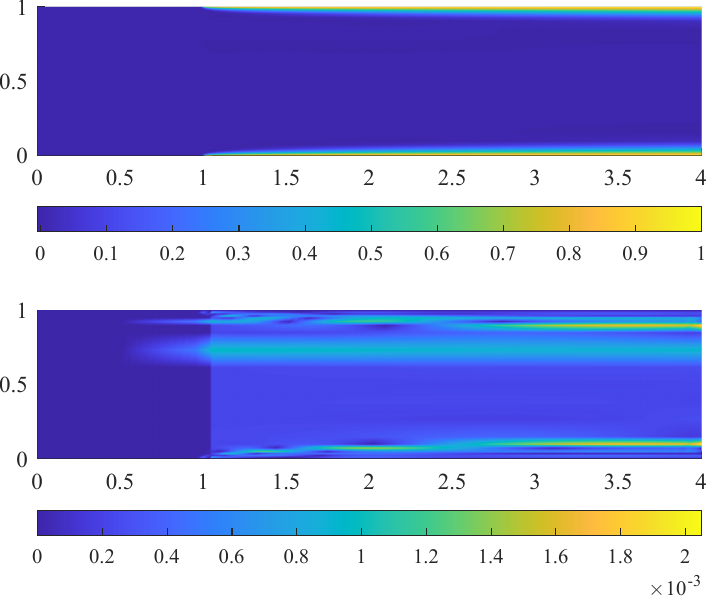}
    \caption{\rev{DD-PGD solution $\upgd(\bmu)$ (top) and the scaled error $| \upgd(\bmu) - u_{\Omega}^h(\bmu) | / \max_{\Omega} |u_{\Omega}^h(\bmu)|$ (bottom) for $\bmu = (1.25 \times 10^4, 3)$.}}
    \label{fig:PGFlow_12500_3}
\end{figure}

\begin{figure}[h!]
    \centering
    \includegraphics[]{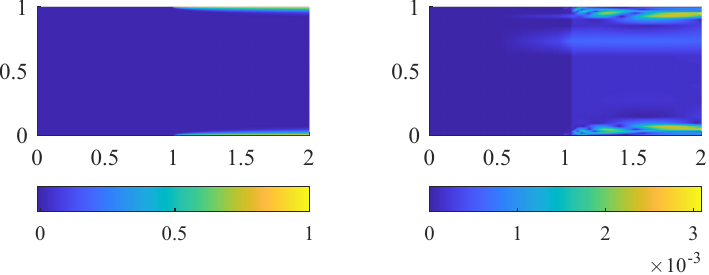}
    \caption{\rev{DD-PGD solution $\upgd(\bmu)$ (left) and the scaled error $| \upgd(\bmu) - u_{\Omega}^h(\bmu) | / \max_{\Omega} |u_{\Omega}^h(\bmu)|$ (right) for $\bmu = (2 \times 10^4, 1)$.}}
    \label{fig:PGFlow_20000_1}
\end{figure}

\rev{
To assess the capability of the DD-PGD strategy to accurately compute the solution of the parametric problem \eqref{eq:pgProblem}, the online evaluation of the surrogate model $\upgd(\bmu)$ is presented, together with the scaled error map $| \upgd(\bmu) - u_{\Omega}^h(\bmu) | / \max_{\Omega} |u_{\Omega}^h(\bmu)|$, $u_{\Omega}^h(\bmu)$ being the corresponding monolithic finite element solution. Two pairs of values $(\mu_1, \mu_2)$ are tested: Fig.~\ref{fig:PGFlow_12500_3} shows the case $\bmu = (1.25 \times 10^4, 3)$, reproducing the results in~\cite{Pacciarini:2014:CMAME}, while Fig.~\ref{fig:PGFlow_20000_1} the case $\bmu = (2 \times 10^4, 1)$. In both cases, the outcome of the DD-PGD algorithm provides a smooth solution, showing excellent agreement with the high-fidelity simulation. Indeed, the maximum value of the scaled error mentioned above achieves $2 \times 10^{-3}$ for $\bmu = (1.25 \times 10^4, 3)$ and $3 \times 10^{-3}$ for $\bmu = (2 \times 10^4, 1)$.
}

Setting the maximum dimension of the Krylov subspace for the GMRES algorithm to 8, the computed DD-PGD solutions converge in 10 and 12 iterations for $\bmu = (1.25 \times 10^4, 3)$ and $\bmu = (2 \times 10^4, 1)$, respectively, with an average computing time of \SI[parse-numbers=false]{1.45 \times 10^{-1}}{s}.
The online scheme accurately glues the local PGD solutions, achieving a maximum difference in $l^\infty$ norm in the overlapping region of the order of $10^{-3}$.
For the sake of comparison, the standard DD-FEM approach, executed for 100 random pairs $(\mu_1, \mu_2)$, requires an average of 14 iterations to converge with an average CPU time of \SI{18.13}{s}.
Hence, the surrogate-based overlapping Schwarz method is approximately 118 times faster than the standard DD-FEM strategy, showcasing that the DD-PGD method provides a competitive framework, also in the presence of convection-dominated phenomena and geometrically parametrized domains.

\subsection{Multi-domain thermal problem with discontinuous conductivity}
\label{sec:testPatera}

\begin{figure}[h!]
    \centering
    \begin{tikzpicture}
        \def\s{2}
        \begin{scope}
            \pateraWing{1}{0}{1}{0}
            \pateraBlock{1}{1}{1}{2}{2}{1}{1}
        \end{scope}

        \begin{scope}[shift={(1.5*\s, 0)}]
            \pateraWing{1}{0}{0}{0}
            \pateraBlock{2}{1}{2}{2}{2}{1}{1}
        \end{scope}

        \begin{scope}[shift={(3*\s, 0)}]
            \pateraWing{1}{0}{0}{1}
            \pateraBlock{3}{1}{2}{1}{2}{1}{1}
            \draw[ultra thick, blue] (0, -0.25*\s) -- (1*\s, -0.25*\s);
        \end{scope}

        \begin{scope}[shift={(0, 1.5*\s)}]
            \pateraWing{0}{0}{1}{0}
            \pateraBlock{4}{2}{1}{2}{2}{1}{1}
        \end{scope}

        \begin{scope}[shift={(1.5*\s, 1.5*\s)}]
            \pateraWing{0}{0}{0}{0}
            \pateraBlock{5}{2}{2}{2}{2}{1}{1}
        \end{scope}

        \begin{scope}[shift={(3*\s, 1.5*\s)}]
            \pateraWing{0}{0}{0}{1}
            \pateraBlock{6}{2}{2}{1}{2}{1}{1}
        \end{scope}

        \begin{scope}[shift={(0, 3*\s)}]
            \pateraWing{0}{1}{1}{0}
            \pateraBlock{7}{2}{1}{2}{1}{1}{1}
            \draw[ultra thick, red] (0, 1.25*\s) -- (1*\s, 1.25*\s);
        \end{scope}

        \begin{scope}[shift={(1.5*\s, 3*\s)}]
            \pateraWing{0}{1}{0}{0}
            \pateraBlock{8}{2}{2}{2}{1}{1}{1}
        \end{scope}

        \begin{scope}[shift={(3*\s, 3*\s)}]
            \pateraWing{0}{1}{0}{1}
            \pateraBlock{9}{2}{2}{1}{1}{1}{1}
        \end{scope}

        \node at (0.45*\s, 4.4*\s) {$\Gamma_{\text{in}}$};

        \node at (3.5*\s, -0.4*\s) {$\Gamma_{\text{out}}$};

    \end{tikzpicture}
    \caption{Modular domain $\Omega$, featuring local bulk regions $\Omega_i^b$ ($i=1,\ldots,9$) with conductivity $\mu_i$ (pink) and wing regions (grey) with unitary conductivity.}
    \label{fig:patera9}
\end{figure}
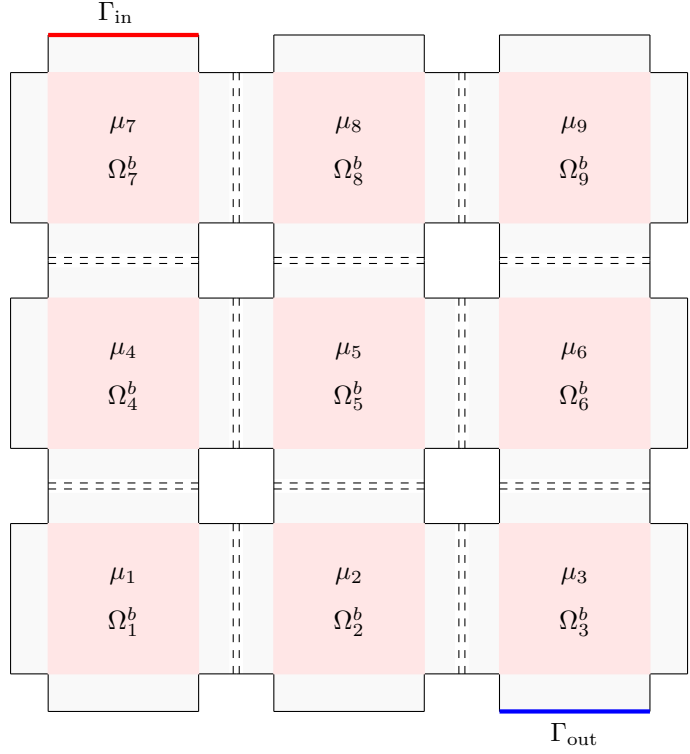

The last test case discusses the benchmark problem introduced in~\cite{Eftang:2013:IJNME} of the parametric thermal equation
\begin{equation}\label{eq:PateraProblem}
    \begin{array}{rcll}
        -\nabla \cdot (\nu(\bmu) \nabla u(\bmu)) &=& 0 & \quad \text{in } \Omega \, ,\\
        u(\bmu) &=& 0 & \quad \text{on } \Gamma_{\text{out}} \, , \\
        \nu(\bmu)\nabla u(\bmu) \cdot \bn &=& 1  & \quad \text{on } \Gamma_{\text{in}} \, ,\\
        \nabla u(\bmu) \cdot \bn &=& 0 & \quad \text{on } \partial\Omega \setminus (\Gamma_{\text{in}} \cup \Gamma_{\text{out}}) \, ,
    \end{array}
\end{equation}
where $\Omega$ is the modular structure consisting of 9 subdomains shown in Fig.~\ref{fig:patera9} and $\nu(\bmu)$ denotes the space-dependent thermal conductivity
\begin{equation}\label{eq:conductivity}
\nu(\bmu) = \begin{cases}
\mu_i & \text{in $\Omega_i^b$ ($i =1,\ldots,9$)}, \\
1 & \text{otherwise} \, .
\end{cases}
\end{equation}

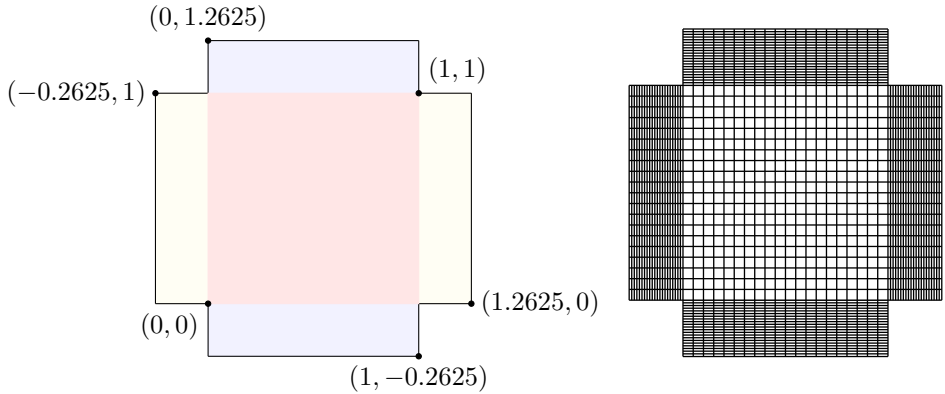
\begin{figure}[h!]
    \centering
    \begin{tikzpicture}
        \def\s{2.8}
        
        \pateraWing{2}{2}{2}{2}
        \pateraBlock{}{1}{1}{1}{1}{0}{0}

        \filldraw[black] (0, 0) circle (1pt) node[anchor=north east]{$(0, 0)$};
        \filldraw[black] (1*\s, 1*\s) circle (1pt) node[anchor=south west]{$(1, 1)$};

        \filldraw[black] (1.25*\s, 0) circle (1pt) node[anchor=west]{$(1.2625, 0)$};
        \filldraw[black] (1*\s, -0.25*\s) circle (1pt) node[anchor=north]{$(1, -0.2625)$};
        \filldraw[black] (-0.25*\s, 1*\s) circle (1pt) node[anchor=east]{$(-0.2625, 1)$};
        \filldraw[black] (0, 1.25*\s) circle (1pt) node[anchor=south]{$(0, 1.2625)$};
    \end{tikzpicture}
    \includegraphics[]{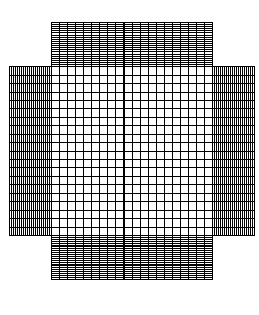}
    \caption{Geometric description of one subdomain (left) and computational mesh (right).}
    \label{fig:pateraMesh}
\end{figure}

The geometric configuration of each subdomain is reported in Fig.~\ref{fig:pateraMesh}. A structured mesh of 2,080 quadrilateral elements is considered, with uniform size $\hX^b = 5 \times 10^{-2}$ in the bulk region (pink), non-uniform mesh size with horizontal spacing $h_{x,1}^v = 5 \times 10^{-2}$ and vertical spacing $h_{x,2}^v = 1.25 \times 10^{-2}$ in the vertical wings (yellow), whereas in the horizontal wings (blue) the spacing in the $x$ and $y$ directions is given by $h_{x,1}^h = 1.25 \times 10^{-2}$ and $h_{x,2}^h = 5 \times 10^{-2}$, respectively.

It is worth noticing that the resulting problem features two spatial dimensions and 9 independent physical parameters, one for each subdomain, leading to a problem of dimension 11. Hence, the construction of a monolithic ROM for equation~\eqref{eq:PateraProblem} is particularly challenging. Nonetheless, as observed in Sect.~\ref{sec:testRozza}, the DD-PGD approach can be utilized to partition the domain into subregions where only a subset of the parameters is relevant. 
Following the rationale originally proposed in~\cite{Eftang:2013:IJNME}, the modular nature of the domain is exploited to partition $\Omega$ into 9 regions, each featuring a problem with a unique physical parameter $\hmu$.

\begin{figure}[h!]
    \centering
    \begin{subfigure}{0.49\textwidth}
        \centering
        \begin{tikzpicture}
            \def\s{2}
            \pateraWing{1}{1}{1}{1}
            \pateraBlock{1}{0}{0}{1}{1}{2}{2}

            \draw[ultra thick, blue] (-0.25*\s, 0) -- (-0.25*\s, 1*\s);
            \node at (-0.7*\s, 0.5*\s) {$u = \lambda_1$};

            \draw[ultra thick, blue] (0, -0.25*\s) -- (1*\s, -0.25*\s);
            \node at (0.5*\s, -0.5*\s) {$u = \lambda_2$};

            \node[white] at (1.75*\s, 0.5*\s) {$u = \lambda_3$};

            \draw[ultra thick, green] (0, 1.25*\s) -- (1*\s, 1.25*\s);
            \draw[ultra thick, green] (1.25*\s, 0) -- (1.25*\s, 1*\s);
            \draw[ultra thick, green] (-0.25*\s, 1*\s) -- (0, 1*\s);
            \draw[ultra thick, green] (1*\s, 1*\s) -- (1.25*\s, 1*\s);
            \draw[ultra thick, green] (-0.25*\s, 0) -- (0, 0);
            \draw[ultra thick, green] (1*\s, 0) -- (1.25*\s, 0);
            \draw[ultra thick, green] (0, -0.25*\s) -- (0, 0);
            \draw[ultra thick, green] (0, 1*\s) -- (0, 1.25*\s);
            \draw[ultra thick, green] (1*\s, -0.25*\s) -- (1*\s, 0);
            \draw[ultra thick, green] (1*\s, 1*\s) -- (1*\s, 1.25*\s);
        \end{tikzpicture}
        \caption{Reference subdomain $\hOmega_1$}
    \end{subfigure}
    \hfill
    \begin{subfigure}{0.49\textwidth}
        \centering
        \begin{tikzpicture}
            \def\s{2}
            \pateraWing{1}{1}{1}{1}
            \pateraBlock{2}{0}{0}{0}{1}{2}{2}

            \draw[ultra thick, blue] (-0.25*\s, 0) -- (-0.25*\s, 1*\s);
            \node at (-0.7*\s, 0.5*\s) {$u = \lambda_1$};

            \draw[ultra thick, blue] (0, -0.25*\s) -- (1*\s, -0.25*\s);
            \node at (0.5*\s, -0.5*\s) {$u = \lambda_2$};

            \draw[ultra thick, blue] (1.25*\s, 0) -- (1.25*\s, 1*\s);
            \node at (1.7*\s, 0.5*\s) {$u = \lambda_3$};

            \draw[ultra thick, green] (0, 1.25*\s) -- (1*\s, 1.25*\s);
            \draw[ultra thick, green] (-0.25*\s, 1*\s) -- (0, 1*\s);
            \draw[ultra thick, green] (1*\s, 1*\s) -- (1.25*\s, 1*\s);
            \draw[ultra thick, green] (-0.25*\s, 0) -- (0, 0);
            \draw[ultra thick, green] (1*\s, 0) -- (1.25*\s, 0);
            \draw[ultra thick, green] (0, -0.25*\s) -- (0, 0);
            \draw[ultra thick, green] (0, 1*\s) -- (0, 1.25*\s);
            \draw[ultra thick, green] (1*\s, -0.25*\s) -- (1*\s, 0);
            \draw[ultra thick, green] (1*\s, 1*\s) -- (1*\s, 1.25*\s);
        \end{tikzpicture}
        \caption{Reference subdomain $\hOmega_2$}
    \end{subfigure}
    \hfill
    \begin{subfigure}{0.49\textwidth}
        \centering
        \begin{tikzpicture}
            \def\s{2}
            \pateraWing{1}{1}{1}{1}
            \pateraBlock{3}{0}{0}{0}{0}{2}{2}

            \draw[ultra thick, blue] (-0.25*\s, 0) -- (-0.25*\s, 1*\s);
            \node at (-0.7*\s, 0.5*\s) {$u = \lambda_1$};

            \draw[ultra thick, blue] (0, -0.25*\s) -- (1*\s, -0.25*\s);
            \node at (0.5*\s, -0.5*\s) {$u = \lambda_2$};

            \draw[ultra thick, blue] (1.25*\s, 0) -- (1.25*\s, 1*\s);
            \node at (1.7*\s, 0.5*\s) {$u = \lambda_3$};

            \draw[ultra thick, blue] (0, 1.25*\s) -- (1*\s, 1.25*\s);
            \node at (0.5*\s, 1.4*\s) {$u = \lambda_4$};

            %
            \draw[ultra thick, green] (-0.25*\s, 1*\s) -- (0, 1*\s);
            \draw[ultra thick, green] (1*\s, 1*\s) -- (1.25*\s, 1*\s);
            \draw[ultra thick, green] (-0.25*\s, 0) -- (0, 0);
            \draw[ultra thick, green] (1*\s, 0) -- (1.25*\s, 0);
            \draw[ultra thick, green] (0, -0.25*\s) -- (0, 0);
            \draw[ultra thick, green] (0, 1*\s) -- (0, 1.25*\s);
            \draw[ultra thick, green] (1*\s, -0.25*\s) -- (1*\s, 0);
            \draw[ultra thick, green] (1*\s, 1*\s) -- (1*\s, 1.25*\s);
        \end{tikzpicture}
        \caption{Reference subdomain $\hOmega_3$}
    \end{subfigure}
    \hfill
    \begin{subfigure}{0.49\textwidth}
        \centering
        \begin{tikzpicture}
            \def\s{2}
            \pateraWing{1}{1}{1}{1}
            \pateraBlock{4}{0}{1}{0}{1}{2}{2}

            \draw[ultra thick, blue] (0, -0.25*\s) -- (1*\s, -0.25*\s);
            \node at (0.5*\s, -0.5*\s) {$u = \lambda_1$};

            \draw[ultra thick, blue] (1.25*\s, 0) -- (1.25*\s, 1*\s);
            \node at (1.7*\s, 0.5*\s) {$u = \lambda_2$};

            \draw[ultra thick, red] (0, 1.25*\s) -- (1*\s, 1.25*\s);
            \node at (0.45*\s, 1.4*\s) {$\nabla u\cdot \bn=1$};

            \node[white] at (-0.7*\s, 0.5*\s) {$u = \lambda_1$};

            \draw[ultra thick, green] (-0.25*\s, 0) -- (-0.25*\s, 1*\s);
            \draw[ultra thick, green] (-0.25*\s, 1*\s) -- (0, 1*\s);
            \draw[ultra thick, green] (1*\s, 1*\s) -- (1.25*\s, 1*\s);
            \draw[ultra thick, green] (-0.25*\s, 0) -- (0, 0);
            \draw[ultra thick, green] (1*\s, 0) -- (1.25*\s, 0);
            \draw[ultra thick, green] (0, -0.25*\s) -- (0, 0);
            \draw[ultra thick, green] (0, 1*\s) -- (0, 1.25*\s);
            \draw[ultra thick, green] (1*\s, -0.25*\s) -- (1*\s, 0);
            \draw[ultra thick, green] (1*\s, 1*\s) -- (1*\s, 1.25*\s);
        \end{tikzpicture}
        \caption{Reference subdomain $\hOmega_4$}
    \end{subfigure}
    \caption{Reference subdomains $\hOmega_j$ ($j=1,\ldots,4$) used to build the modular domain $\Omega$ and corresponding boundary types: Dirichlet (blue), unitary Neumann (red) and homogeneous Neumann (green) boundary conditions.}
    \label{fig:patera9Ref}
\end{figure}
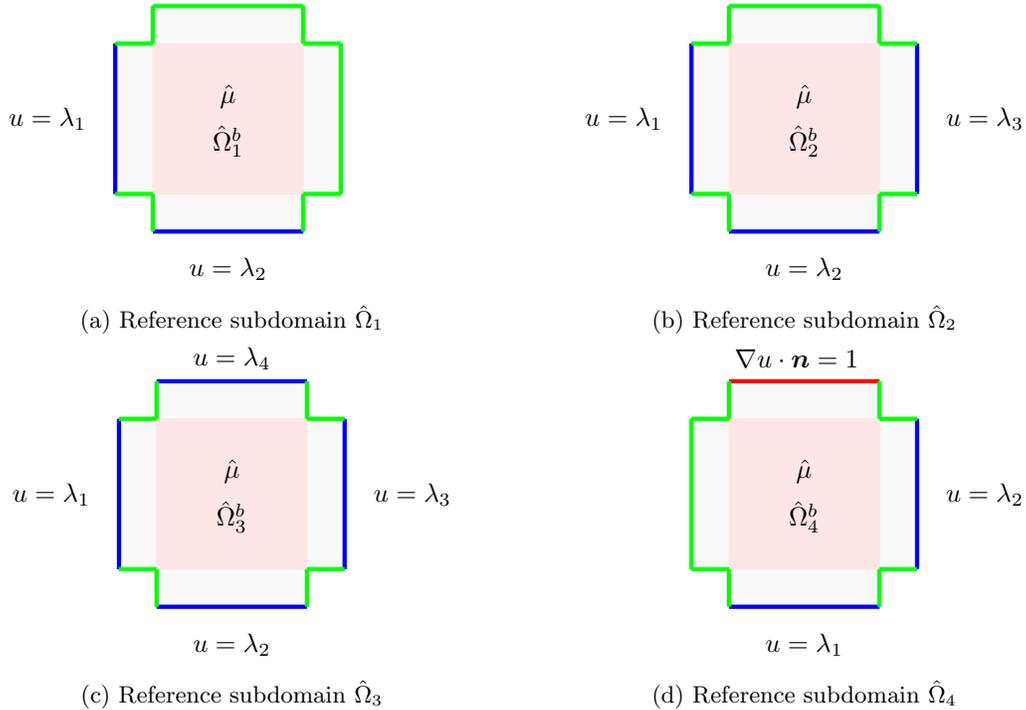

The domain $\Omega$ is thus split into 9 overlapping subdomains $\Omega_i$ ($i=1,\ldots,9$) with overlap width equal to $2.5 \times 10^{-2}$. \rev{Note that the overlapping regions are located in the horizontal and vertical wings of the subdomains, where the thermal conductivity is maintained constant and equal to 1, avoiding any possible issue related to overlaps with discontinuous material properties.} Considering the type of conditions imposed on the boundary of each subdomain, four reference subdomains $\hOmega_j$ ($j=1,\ldots,4$) can be identified that, upon appropriate translation and/or rotation, can be used to describe all subdomains $\Omega_i$ ($i=1,\ldots,9$) thus reconstructing the original domain $\Omega$. The reference subdomains are shown in Fig.~\ref{fig:patera9Ref} with the corresponding imposed boundary conditions. When no condition is specified, a homogeneous Neumann boundary condition is applied. Table~\ref{tab:onlineTransformations} reports the transformations that must be applied to $\hOmega_j$ ($j=1,\ldots,4$) to retrieve $\Omega_i$ ($i=1,\ldots,9$).

\begin{table}[h!]
\begin{center}
\resizebox{\textwidth}{!}{%
\begin{tabular}{lccccccccc}
            & $\Omega_1$  & $\Omega_2$  & $\Omega_3$  & $\Omega_4$  & $\Omega_5$  & $\Omega_6$  & $\Omega_7$  & $\Omega_8$  & $\Omega_9$ \\
\hline \\[-1em]
Ref. subdomain & $\hOmega_1$ & $\hOmega_2$ & $\hOmega_2$ & $\hOmega_2$ & $\hOmega_3$ & $\hOmega_2$ & $\hOmega_4$ & $\hOmega_2$ & $\hOmega_1$ \\[0.25em]
Translation & $\begin{pmatrix}0 \\ 0\end{pmatrix}$ & $\begin{pmatrix} 1.5 \\ 0\end{pmatrix}$ & $\begin{pmatrix}3 \\ 0\end{pmatrix}$ & $\begin{pmatrix}0 \\ 1.5 \end{pmatrix}$ & $\begin{pmatrix} 1.5 \\ 1.5 \end{pmatrix}$ & $\begin{pmatrix} 3 \\ 1.5\end{pmatrix}$ & $\begin{pmatrix} 0 \\ 3\end{pmatrix}$ & $\begin{pmatrix} 1.5 \\ 3 \end{pmatrix}$ & $\begin{pmatrix} 3 \\ 3 \end{pmatrix}$ \\[1em]
Rotation    & $\pi$ & $\pi$ & $\displaystyle\frac{3\pi}{2}$ & $\displaystyle\frac{\pi}{2}$ & 0 & $\displaystyle\frac{3\pi}{2}$ & 0 & 0 & 0\\[0.5em]
Case study \#1: $\mu_i=$ & 0.1 & 0.2 & 0.4 & 0.8 & 1.6 & 3.2 & 6.4 & 0.1 & 0.2 \\
Case study \#2: $\mu_i=$ & 4.9 & 4.7 & 4.8 & 5.2 & 5 & 4.9 & 5.5 & 5.3 & 5.1 \\
\hline
\end{tabular}}
\end{center}
\caption{Transformation of the reference subdomains $\hOmega_j$ ($j=1,\ldots,4$) to obtain the physical subdomains $\Omega_i$ ($i=1,\ldots,9$), and conductivity parameters $\mu_i$ in the physical subdomains for two case studies.}
\label{tab:onlineTransformations}
\end{table}

The physical parameter $\mu_i$ is set to vary in the interval $\mathcal{P} = [5 \times 10^{-2}, 10]$, while all interface boundary parameters belong to the interval $\mathcal{J}_i^q = [-5, 5]$, for all $i$ and $q$. The parametric intervals are discretized using equally-space nodes with mesh size $h_{\mu} = \hLam = 10^{-3}$, yielding a total of $N^{\mu} = 9.95 \times 10^3$ and $\Nlam = 10^4$ unknowns in $\mathcal{P}$ and $\mathcal{J}_i^q$, respectively.

The spatial problem is discretized using continuous $\mathbb{Q}_2$ Lagrangian finite element functions. For each subdomain, Table~\ref{tab:Patera} reports the resulting number of finite element degrees of freedom after taking into account the boundary conditions, the number of physical and boundary parameters and the total number of PGD modes computed for the local surrogate models. 
Due to the presence of non-homogeneous boundary data only on $\Gamma_{\text{in}}$, the problem of type~\eqref{eq:sourceProb} needs to be formulated solely in the reference subdomain $\hOmega_4$. Concerning the problems of type~\eqref{eq:boundaryProbs}, the active boundary parameters are partitioned to ensure that each subproblem features at most three parameters, yielding 14 subproblems in the reference subdomain $\hOmega_1$, 21 in $\hOmega_2$, 28 in $\hOmega_3$ and 14 in $\hOmega_4$. 
The overall computational time of the offline phase is \SI{630}{s}, with additional \SI{228}{s} required to perform the PGD compression.
\begin{table}[bht]
\begin{center}
\begin{tabular}{cccccc}
Reference & Finite & Physical & Interface boundary & PGD & Modes after\\
subdomain & element DOFs & parameters & parameters\footnotemark & modes & compression\\
\hline
$\hOmega_1$ & 2,163 & 1 & 42 & 370 & 236\\
$\hOmega_2$ & 2,142 & 1 & 63 & 672 & 386\\
$\hOmega_3$ & 2,121 & 1 & 84 & 930 & 509\\
$\hOmega_4$ & 2,163 & 1 & 42 & 419 & 244\\
\hline
\end{tabular}
\end{center}
\caption{Specifics of the local subproblems and of the computed local PGD surrogate models for each reference subdomain.}
\label{tab:Patera}
\end{table}
\footnotetext{When the reference subdomain $\hOmega_2$ is employed for the computation of the physical subdomain $\Omega_3$, only 42 interface boundary parameters are considered, the remaining 21 being fixed by imposing the homogeneous Dirichlet condition on $\Gamma_{\text{out}}$.}

In the online phase, the solution in the global domain $\Omega$ is reconstructed by composing the local PGD surrogate models computed in the reference subdomains $\hOmega_j$ ($j=1,\ldots,4$) with the appropriate translations and rotations described in Table~\ref{tab:onlineTransformations}. In total, there are 12 overlapping regions between the nine subdomains, identified by the dashed lines in Fig.~\ref{fig:patera9}, and the dimension of the resulting interface system is 504.

Two case studies associated with the different values of the conductivity parameters reported in Table~\ref{tab:onlineTransformations} are considered.
Case study \#1 reproduces the benchmark in~\cite{Eftang:2013:IJNME}.
Setting the maximal dimension of the Krylov subspace equal to 60, GMRES performs 302 iterations, converging in approximately \SI{19.13}{s}. The corresponding DD-FEM solver requires 95 GMRES iterations to converge, with a total CPU time of \SI{414.8}{s}. As already observed in the previous numerical experiments, despite the larger number of GMRES iterations, the DD-PGD method still outperforms the standard DD-FEM procedure in terms of computing time, with an approximate speed-up of 22 times.

\begin{figure}[h!]
    \centering
    \includegraphics[]{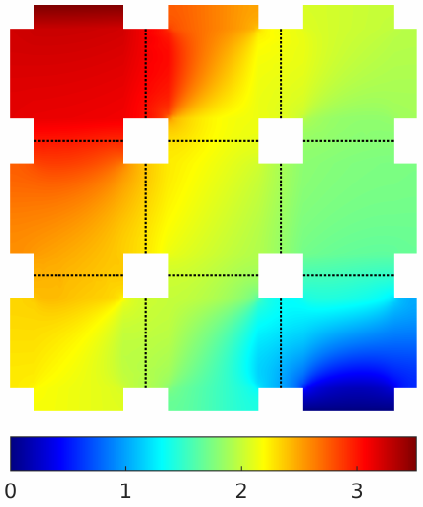}
    \hfill
    \includegraphics[]{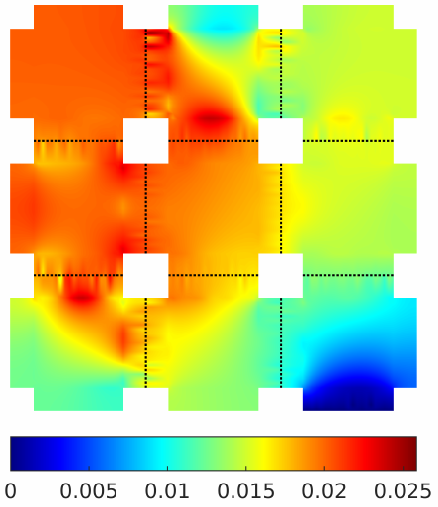}
    \caption{\rev{Global DD-PGD solution (left) and the scaled error $| \upgd(\bmu) - u^h(\bmu) | / \max_{\Omega} |u^h(\bmu)|$ (right) for case study \#1 in Table~\ref{tab:onlineTransformations}. The dashed lines indicate the location of the overlapping regions.}}
    \label{fig:patera9comparison}
\end{figure}

Figure~\ref{fig:patera9comparison} displays the temperature distribution $\upgd(\bmu)$ obtained using the DD-PGD strategy and \rev{the scaled  error $| \upgd(\bmu) - u^h(\bmu) | / \max_{\Omega} |u^h(\bmu)|$ for the case study \#1. The maximum value of the scaled error is $2.5 \times 10^{-2}$.}
Moreover, the difference between the local surrogate models in the overlapping region, measured in the $l^\infty$ norm, is equal to $5 \times 10^{-9}$, showcasing the accurate imposition of the continuity of solution in the online phase of the DD-PGD method, even in the case of partitions involving multiple subdomains. 

\begin{figure}[h!]
    \centering
    \begin{subfigure}[][80.5mm][b]{0.49\textwidth}
        \centering
        \includegraphics{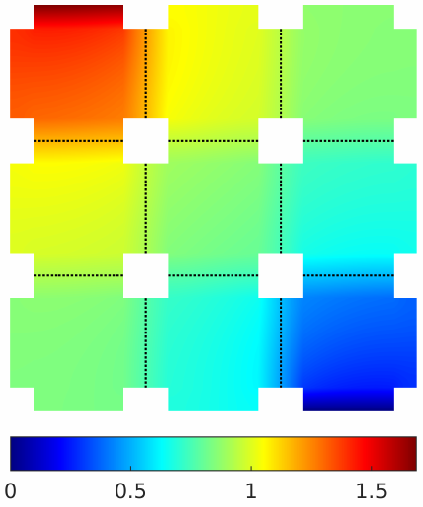}
    \end{subfigure}
    \hfill
    \begin{subfigure}[][][t]{0.49\textwidth}
        \centering
        \includegraphics{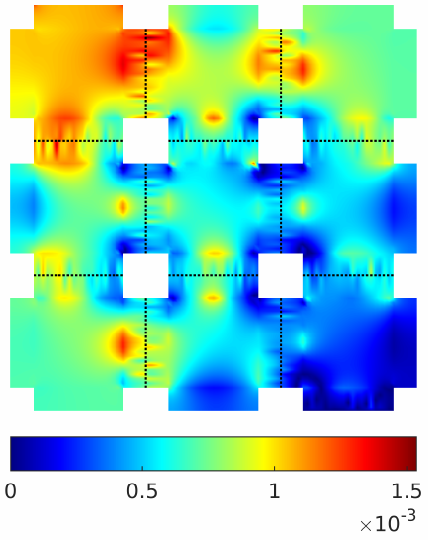}
    \end{subfigure}
    \caption{\rev{Global DD-PGD solution (left) and the scaled error $| \upgd(\bmu) - u^h(\bmu) | / \max_{\Omega} |u^h(\bmu)|$ (right) for case study \#2 in Table~\ref{tab:onlineTransformations}. The dashed lines indicate the location of the overlapping regions.}}
    \label{fig:pateraAlt9comparison}
\end{figure}

A second test is carried out considering the alternative choice of conductivity parameters identified as case study \#2 in Table~\ref{tab:onlineTransformations}. This choice avoids large variations of the conductivity between adjacent subdomains and guarantees that the diffusion operator is far from becoming singular.
In this case, the convergence behaviour of GMRES significantly improves for both the DD-PGD and the DD-FEM approaches. More precisely, for DD-PGD, GMRES converges in 126 iterations with an average CPU time of \SI{8.32}{s}, while, for DD-FEM, convergence is achieved in \SI{249}{s} and 56 iterations. Hence, for case study \#2, the DD-PGD method is approximately 30 times faster than the standard DD-FEM procedure. The comparison of the DD-PGD and DD-FEM temperature distributions for case study \#2 is presented in Fig.~\ref{fig:pateraAlt9comparison}, showing the excellent accuracy of the ROM solution \rev{for which the maximum value of the scaled error $| \upgd(\bmu) - u^h(\bmu) | / \max_{\Omega} |u^h(\bmu)|$ is $1.5 \times 10^{-3}$.}

The improvement in the convergence behaviour observed in case study \#2 suggests that the robustness of the iterative method can be enhanced by introducing \emph{ad-hoc} preconditioning strategies to tackle both the variations in the physical parameters and the presence of multiple subdomains. This optimization of the online phase, which is beyond the contributions of this paper, could be achieved following several approaches from domain decomposition, see, e.g., \cite{Toselli:2005,Dolean-DJN-15} and references therein.

\section{Concluding remarks}
\label{sec:Conclusions}

In this work, a DD-PGD approach combining the overlapping Schwarz algorithm with a physics-based PGD-ROM was proposed to solve parametric linear elliptic PDEs.
The method constructs local surrogate models with low dimensionality by exploiting the linearity of the problems to define disjoint sets of active boundary parameters that allow to represent arbitrary Dirichlet boundary conditions at the interfaces between subdomains.
In the online phase, the coupling of the subdomains is performed via an interface equation, leading to a linear system for the nodal values of the parametric solution at the interfaces.
The resulting system can be efficiently solved in real time using standard Krylov methods.

The advantages of the proposed approach are multiple.
First, by relying on a fully algebraic formulation, the DD-PGD strategy provides a non-intrusive framework for the construction of local ROMs via the encapsulated PGD library.
In the offline phase, the traces of the finite element functions used for the spatial discretization within each subdomain are employed to define parametric Dirichlet boundary conditions at the subdomain level, without the need to introduce auxiliary basis functions for the solution at the interface. 
In addition, the linearity of the operators is exploited to reduce the overall dimensionality of the problem, by devising a set of subproblems with only a few active boundary parameters.
For the coupling procedure, the discussed parametric multi-domain formulation allows to seamlessly glue the local ROMs only at the interfaces, without introducing extra variables (i.e., Lagrange multipliers) or enforcing the continuity of the local solutions in the entire overlapping region.
Finally, the solution of the parametric interface equation only entails the interpolation of the previously computed surrogate models in the parametric space, with no additional problems to be solved during the online phase.

The resulting DD-PGD approach was tested on a set of numerical benchmarks, including one and multiple parameters (both physical and geometrical), two and multiple subdomains, to assess accuracy, robustness and efficiency of the method.
The strategy showed accurate results comparable with the high-fidelity solutions, and robustness in different scenarios (from pure diffusion to convection-dominated convection-diffusion equations) while outperforming the standard non-overlapping DD-FEM in terms of computing time.
\rev{
It is worth noticing that the CPU times reported for the DD-PGD simulations could be further optimized by accelerating the convergence of the online phase via tailored preconditioning strategies~\cite{Toselli:2005,Dolean-DJN-15}. Other robust overlapping DD algorithms as well as non-overlapping DD strategies, which are outside the scope of this work, should also be investigated. These aspects play a crucial role in guaranteeing the applicability of the described methodology to real-world, three-dimensional cases. Indeed, although the presented approach can be seamlessly employed to construct local surrogate models of 3D problems, the partition of complex domains of industrial interest relies on state-of-the-art graph partitioning software, such as, e.g., METIS, KaHIP and Scotch~\cite{Metis,Kahip,Scotch}, and their efficient, non-intrusive coupling with local PGD-ROMs will require further study.
}

\paragraph{Acknowledgements} The authors acknowledge funding as follows. MD: EPSRC grant EP/V027603/1. BJE: EPSRC Doctoral Training Partnership grant EP/W523987/1. MG: Spanish Ministry of Science and Innovation and Spanish State Research Agency \\ MCIN/AEI/10.13039/501100011033 (Grants No.  PID2020-113463RB-C33 and CEX2018-000797-S).  MG is Fellow of the Serra H\'unter Programme of the Generalitat de Catalunya.

\bibliographystyle{elsarticle-num}
\bibliography{references_paper_ddpgd_elliptic}

\appendix

\section{Encapsulated proper generalized decomposition solver}\label{sec:encapsulated}
\label{append:encapsulatedPGD}

In this appendix, some technical aspects related to the solution of the parametric linear systems~\eqref{eq:algebraicA} and~\eqref{eq:algebraicB} using the encapsulated PGD toolbox are presented. For a detailed description of the method, the interested reader is referred to~\cite{Diez:2020:ACME}.

The goal is to compute the solutions $\bvpgd_{i,0}(\bmu)$ and $\bvpgd_{i,j}(\bmu, \bLambda_i^j)$ in separated form, that is, as the sum of rank-one approximations.
It follows that $M_0$ pairs $(\bV_{i,0}^m,\phi_{i,0}^m(\bmu))$ need to be determined for the PGD expansion~\eqref{eq:algebraicSolA} of $\bvpgd_{i,0}(\bmu)$. Similarly, the triplets $(\bV_{i,j}^m,\phi_{i,j}^m(\bmu),\psi_{i,j}^m(\bLambda_i^j))$, for $m=1,\ldots,M_j$, are sought to approximate $\bvpgd_{i,j}(\bmu, \bLambda_i^j)$ according to~\eqref{eq:algebraicSolB}.
To this end, the encapsulated PGD method relies on a greedy procedure to compute the $m$-th spatial and parametric modes, assuming that the terms up to $m-1$ are known, thus leading to the parametric problems
\begin{subequations}\label{eq:algebraicGreedy}
\begin{align}
\left(
\sum_{\ell=1}^{n_{\nu}} \xi_{\nu}^\ell(\bmu) \mat{K}_{\nu}^\ell
+ \sum_{\ell=1}^{n_{\alpha}} \xi_{\alpha}^\ell(\bmu) \mat{K}_{\alpha}^\ell
+ \sum_{\ell=1}^{n_{\gamma}} \xi_{\gamma}^\ell(\bmu) \mat{K}_{\gamma}^\ell
\right)
\bV_{i,0}^m \phi_{i,0}^m(\bmu) 
=
\mathbf{R}_{i,0}^{m-1}(\bmu) &
\notag \\[-1em]
\qquad \forall \bmu \in \mathcal{P} \, , &
\label{eq:greedyA} \\
\left(
\sum_{\ell=1}^{n_{\nu}} \xi_{\nu}^\ell(\bmu) \mat{K}_{\nu}^\ell
+ \sum_{\ell=1}^{n_{\alpha}} \xi_{\alpha}^\ell(\bmu) \mat{K}_{\alpha}^\ell
+ \sum_{\ell=1}^{n_{\gamma}} \xi_{\gamma}^\ell(\bmu) \mat{K}_{\gamma}^\ell
\right) 
\bV_{i,j}^m \phi_{i,j}^m(\bmu) \psi_{i,j}^m(\bLambda_i^j ) 
=
\mathbf{R}_{i,j}^{m-1}(\bmu, \bLambda_i^j ) &
\notag \\[-1em]
\forall \bmu \in \mathcal{P} \, \text{ and } \forall \bLambda_i^j \in \mathcal{Q}_i^j \, , &
\label{eq:greedyB}
\end{align}
\end{subequations}
where $\mathbf{R}_{i,0}^{m-1}(\bmu)$ and $\mathbf{R}_{i,j}^{m-1}(\bmu, \bLambda_i^j )$ denote the residuals
\begin{subequations}\label{eq:residualsGreedy}
\begin{align}
\mathbf{R}_{i,0}^{m-1}(\bmu) :=&
\sum_{\ell=1}^{n_s} \xi_s^\ell(\bmu) \mathbf{f}_s^\ell
+ \sum_{\ell=1}^{n_N} \xi_N^\ell(\bmu) \mathbf{f}_N^\ell
+ \sum_{\ell=1}^{n_D} \xi_D^\ell(\bmu) \mathbf{f}_D^\ell
\notag \\
&-
\sum_{k=1}^{m-1} 
\left(
\sum_{\ell=1}^{n_{\nu}} \xi_{\nu}^\ell(\bmu) \mat{K}_{\nu}^\ell
+ \sum_{\ell=1}^{n_{\alpha}} \xi_{\alpha}^\ell(\bmu) \mat{K}_{\alpha}^\ell
+ \sum_{\ell=1}^{n_{\gamma}} \xi_{\gamma}^\ell(\bmu) \mat{K}_{\gamma}^\ell
\right)
\bV_{i,0}^k \phi_{i,0}^k(\bmu) \, ,
\label{eq:resGreedyA} \\
\mathbf{R}_{i,j}^{m-1}(\bmu, \bLambda_i^j ) :=&
\sum_{q \in \mathcal{N}_i^j} \Lambda^q_i \mathbf{f}_{\Lambda}^q 
\notag \\
&-
\sum_{k=1}^{m-1} 
\left(
\sum_{\ell=1}^{n_{\nu}} \xi_{\nu}^\ell(\bmu) \mat{K}_{\nu}^\ell
+ \sum_{\ell=1}^{n_{\alpha}} \xi_{\alpha}^\ell(\bmu) \mat{K}_{\alpha}^\ell
+ \sum_{\ell=1}^{n_{\gamma}} \xi_{\gamma}^\ell(\bmu) \mat{K}_{\gamma}^\ell
\right) 
\bV_{i,j}^k \phi_{i,j}^k(\bmu) \psi_{i,j}^k(\bLambda_i^j ) \, .
\label{eq:resGreedyB}
\end{align}
\end{subequations}

In order to make the solution of high-dimensional equations~\eqref{eq:algebraicGreedy} feasible, the encapsulated PGD framework uses a fixed-point iteration scheme, namely the alternating directions algorithm, whose non-intrusive implementation is described in~\cite{Diez:2020:ACME}.
In particular, this library employs a purely algebraic formulation of the PGD, relying on the so-called \texttt{separatedTensor} structure, where all information is stored in the form of either matrices or vectors, as detailed in Appendix~\ref{append:Storage}.
An example of the setup and computation performed by the encapsulated PGD solver is presented in Appendix~\ref{append:Setup} for the test case in Sect.~\ref{sec:testAnalytical}.

\subsection{Encapsulated PGD: solution and data storage}
\label{append:Storage}		

Consider a finite element discretization for the spatial subproblems and a pointwise collocation approach for each parameter. For the sake of simplicity, parameters $\mu$ and $\Lambda_i^j$ are assumed to be scalar but this can be straightforwardly generalized to the vectorial case according to~\eqref{eq:vectParam}.
More precisely, let $\Nfem$ be the number of unknowns arising from the finite element discretization of the spatial variables defined on $\Omega_i$. Similarly,  $\Nmu$ and $\Nlam$ denote the number of unknowns in the uniform meshes introduced for the parametric domains $\mathcal{P}$ and $\mathcal{Q}_i^j$, respectively.

As previously explained, the computation of the $m$-th mode of the PGD expansion of $\bvpgd_{i,0}$ requires determining the pair $(\bV_{i,0}^m,\phi_{i,0}^m)$, depending on $\bx$ and $\mu$, respectively. Similarly, the triplet $(\bV_{i,j}^m,\phi_{i,j}^m,\psi_{i,j}^m)$ of functions of $\bx$, $\mu$ or $\Lambda_i^j$ needs to be computed for $\bvpgd_{i,j}$.
This leads to the spatial modes of the solution being stored in the vectors $\bV_{i,0}^m, \bV_{i,j}^m \in \R^{\Nfem}$ of finite element nodal unknowns in $\Omega_i$, whereas vectors $\bphi_{i,0}^m, \bphi_{i,j}^m \in \R^{\Nmu}$ and $\bpsi_{i,j}^m \in \R^{\Nlam}$ contain the parametric modes of the solution defined for all $\mu \in \mathcal{P}$ and for all $\Lambda_i^j \in \mathcal{Q}_i^j$, respectively.

Similarly, the parametric modes of problem data (i.e.,   $\xi_{\nu}^\ell$, $\xi_{\alpha}^\ell$, $\xi_{\gamma}^\ell$, $\xi_s^\ell$, $\xi_N^\ell$, $\xi_D^\ell$) are discretized for all $\mu \in \mathcal{P}$ and the resulting values are stored in vectors $\bxi_{\nu}^\ell$, $\bxi_{\alpha}^\ell$, $\bxi_{\gamma}^\ell$, $\bxi_s^\ell$, $\bxi_N^\ell$, $\bxi_D^\ell$ of dimension $\Nmu$.
Following the same rationale, for all values of $\Lambda_i^j \in \mathcal{Q}_i^j$, the parametric modes of the active boundary parameters are stored in a set of vectors $\bups_{\Lambda}^1, \bups_{\Lambda}^2, \ldots \in \R^{\Nlam}$.

\subsection{Encapsulated PGD: setup of the test in Section~\ref{sec:testAnalytical}}
\label{append:Setup}

\begin{table}[!ht]
\begin{tabularx}{\textwidth}{lcX}
\hline
\hline\\[-10pt]
Terms in the   & Algebraic   & \\
weak form      & counterpart & 
\\[2pt]
\hline
\multicolumn{3}{c}{Spatial modes} \\
\hline\\[-10pt]
$\displaystyle \int_{\Omega_i} \nabla v \cdot \nabla \deV \,d\bx$ & $\mat{K}_{\nu}^1$ & $\Nfem \times \Nfem$ finite element stiffness matrix with constant diffusion coefficient $b_{\nu}^1 = 1$
\\
$\displaystyle \int_{\Omega_i} x\,\nabla v \cdot \nabla \deV\,d\bx$ & $\mat{K}_{\nu}^2$ & $\Nfem \times \Nfem$ finite element stiffness matrix with space-dependent diffusion coefficient $b_{\nu}^2 = x$
\\
$\displaystyle \int_{\Omega_i} b_s^1(\bx) \deV\,d\bx$ & $\mathbf{f}_s^1$ & $\Nfem \times 1$ vector arising from finite element discretization of the source term $b_s^1$
\\
$\displaystyle \int_{\Omega_i} b_s^2(\bx) \deV\,d\bx$ & $\mathbf{f}_s^2$ & $\Nfem \times 1$ vector arising from finite element discretization of the source term $b_s^2$
\\
$\displaystyle \int_{\Omega_i} b_s^3(\bx) \deV\,d\bx$ & $\mathbf{f}_s^3$ & $\Nfem \times 1$ vector arising from finite element discretization of the source term $b_s^3$
\\
$\displaystyle \int_{\Omega_i} \nabla \varphi^q \cdot \nabla \deV \, d\bx$ & $\mathbf{f}_{\Lambda}^{q,1}$ & $\Nfem \times 1$ vector arising from finite element imposition of Dirichlet boundary conditions
\\[10pt]
$\displaystyle \int_{\Omega_i} x\,\nabla \varphi^q \cdot \nabla \deV \, d\bx$ & $\mathbf{f}_{\Lambda}^{q,2}$ & $\Nfem \times 1$ vector arising from finite element  imposition of Dirichlet boundary conditions
\\[20pt]
\hline
\multicolumn{3}{c}{Parametric modes in $\mu$} \\
\hline\\[-10pt]
$\displaystyle 1$ & $\bxi_{\nu}^1, \,  \bxi_{s}^1$ & $\Nmu \times 1$ vector of ones
\\[10pt]
$\displaystyle \mu$ & $\bxi_{\nu}^2, \,  \bxi_{s}^2$ & $\Nmu \times 1$ vector with discrete values of $\mu \in (\muMin,\muMax)$
\\[10pt]
$\displaystyle \mu^2$ & $\bxi_{s}^3$ & $\Nmu \times 1$ vector with discrete values of $\mu^2$ with $\mu \in (\muMin,\muMax)$
\\[20pt]
\hline
\multicolumn{3}{c}{Parametric modes in $\Lambda$} \\
\hline\\[-10pt]
$1$ & $\bups_{\Lambda}^1$ & $\Nlam \times 1$ vector of ones
\\[10pt]
$\Lambda_i^q$ & $\bups_{\Lambda}^2$ & $\Nlam \times 1$ vector with discrete values of $\Lambda_i^q \in (\LamMin,\LamMax)$
\\[5pt]
\hline
\hline
\end{tabularx}
\caption{Differential and algebraic formulation of the terms in the parametric problem for the setup of the encapuslated PGD solver.}
\label{table:Implementation}
\end{table}

In this appendix, the construction of the \texttt{separatedTensor} structure required by the encapsulated PGD solver is detailed for the test case presented in Sect.~\ref{sec:testAnalytical}.
In this context, a diffusion problem with homogeneous Dirichlet boundary conditions is considered. It follows that $\balpha = \bm{0}$, $\gamma = 0$,  $g^D = 0$, while the Neumann datum $g^N$ is redundant.
From~\eqref{eq:separatedData}, the separated data features $n_{\nu} = 2$ and $n_s=3$ terms, with the spatial and parametric modes given by
\begin{equation}\label{eq:sepDataEx}
\begin{aligned}
\xi_{\nu}^1(\mu) &= 1
\, , \quad
&& b_{\nu}^1(\bx) = 1
\, , \\
\xi_{\nu}^2(\mu) &= \mu
\, , \quad
&& b_{\nu}^2(\bx) = x
\, , \\
\xi_{s}^1(\mu) &= 1
\, , \quad
&& b_{s}^1(\bx) = 8\pi^2\sin(2\pi x)\sin(2\pi y)
\, , \\
\xi_{s}^2(\mu) &= \mu
\, , \quad
&& b_{s}^2(\bx) = 2\pi(4\pi x\sin(2\pi x) - \cos(2\pi x))\sin(2\pi y) - x(x-2) -y(y-1)
\, , \\
\xi_{s}^3(\mu) &= \mu^2
\, , \quad
&& b_{s}^3(\bx) = y(y-1)(1-2x)-x^2(x-2)
\, .
\end{aligned}
\end{equation}

The resulting bilinear form for problems~\eqref{eq:pbApgd} and~\eqref{eq:pbBpgd} is
\begin{equation}\label{eq:bilinearApgdEx}
\Apgd(v, \deV; \mu) 
= \int_{\Omega_i}   \nabla v\cdot \nabla \deV\,d\bx 
+ \mu \int_{\Omega_i}  x \nabla v\cdot \nabla \deV\,d\bx \, ,
\end{equation}
while the linear form for problem~\eqref{eq:pbApgd} is
\begin{subequations}\label{eq:linearPPGDex}
\begin{equation}\label{eq:linearApgdEx}
\Fpgd_0(\deV; \mu)
=  \int_{\Omega_i} b_s^1(\bx) \deV\,d\bx 
+ \displaystyle \mu \int_{\Omega_i} b_s^2(\bx) \deV\,d\bx 
+ \displaystyle \mu^2 \int_{\Omega_i} b_s^3(\bx) \deV\,d\bx
\, ,
\end{equation}
for any value of $\mu \in \mathcal{P}$, and the linear form for problem~\eqref{eq:pbBpgd} is
\begin{equation}\label{eq:linearBpgdEx}
\Fpgd_j(\deV; \mu;  \bLambda_i^j ) = 
- \sum_{q \in \mathcal{N}_i^j} \Lambda^q_i  \left( \int_{\Omega_i}   \nabla \varphi^q_i \cdot \nabla \deV\,d\bx 
+ \mu \int_{\Omega_i}  x \nabla \varphi^q_i\cdot \nabla \deV\,d\bx \right) \,  ,
\end{equation}
for all $\mu \in \mathcal{P}$ and for all $\bLambda_i^j \in \mathcal{Q}_i^j$.
\end{subequations}

It is worth noticing that the integrals in equation~\eqref{eq:bilinearApgdEx} correspond to the standard finite element matrices for the Poisson equation, the first one with a constant unitary diffusion coefficient and the second one with a space-dependent diffusion equal to $x$. Similarly,  the integrals in~\eqref{eq:linearPPGDex} yield the standard finite element vectors on the right-hand side of the linear system, accounting for the source term and the imposition of the Dirichlet boundary conditions.
Table~\ref{table:Implementation} reports a summary of the matrices and vectors required for the construction of the spatial and parametric terms of the \texttt{separatedTensor} structure employed by the encapsulated PGD solver. 

Finally, the setup of the encapsulated PGD solver for problems~\eqref{eq:algebraicGreedy} is briefly presented.
Algorithm~\ref{alg:setupA} details the implementation for problem~\eqref{eq:greedyA}, whereas the solver for problem~\eqref{eq:greedyB} is presented in Algorithm~\ref{alg:setupB}. It is worth recalling that the latter problem features the parametric description of the subdomain boundary conditions, with $\mbox{card}(\mathcal{N}_i^j)$ denoting the number of active boundary parameters.
\begin{algorithm}[!ht]
\caption{Encapsulated PGD solver for problem~\eqref{eq:greedyA}.}\label{alg:setupA}
\begin{algorithmic}[1]
\REQUIRE{Spatial ($\mat{K}_{\nu}^1, \mat{K}_{\nu}^2, \mathbf{f}_s^1, \mathbf{f}_s^2, \mathbf{f}_s^3$) and parametric ($\bxi_{\nu}^1, \bxi_{\nu}^2, \bxi_{s}^1, \bxi_{s}^2, \bxi_{s}^3$) modes of the linear system, tolerances $\varepsilon$ for the PGD enrichment and $\varepsilon^\star$ for the PGD compression.}
\STATE{Initialization of the separated tensor for the left-hand side of problem: \\
{\small \verb|Ki0 = separatedTensor;|}
}
\STATE{Setup of the spatial modes: \\
{\small \verb|Ki0.sectionalData{1,1} =|} $\mat{K}_{\nu}^1${\small \verb|;|} \\
{\small \verb|Ki0.sectionalData{1,2} =|} $\mat{K}_{\nu}^2${\small \verb|;|}
}
\STATE{Setup of the parametric modes: \\
{\small \verb|Ki0.sectionalData{2,1} =|} $\bxi_{\nu}^1${\small \verb|;|} \\
{\small \verb|Ki0.sectionalData{2,2} =|} $\bxi_{\nu}^2${\small \verb|;|} \\
}
\STATE{Initialization of the separated tensor for the right-hand side of the problem: \\
{\small \verb|Fi0 = separatedTensor;|}
}
\STATE{Setup of the spatial modes: \\
{\small \verb|Fi0.sectionalData{1,1} =|} $\mathbf{f}_s^1${\small \verb|;|} \\
{\small \verb|Fi0.sectionalData{1,2} =|} $\mathbf{f}_s^2${\small \verb|;|} \\
{\small \verb|Fi0.sectionalData{1,3} =|} $\mathbf{f}_s^3${\small \verb|;|} \\
}
\STATE{Setup of the parametric modes: \\
{\small \verb|Fi0.sectionalData{2,1} =|} $\bxi_{s}^1${\small \verb|;|} \\
{\small \verb|Fi0.sectionalData{2,2} =|} $\bxi_{s}^2${\small \verb|;|} \\
{\small \verb|Fi0.sectionalData{2,3} =|} $\bxi_{s}^3${\small \verb|;|} \\
}
\STATE{Solution of the problem via the alternating directions method: \\
{\small \verb|vi0 = pgdLinearSolve(Ki0,Fi0,`tolModes',|}$\varepsilon${\small \verb|});|}
}
\IF{PGD compression is active}
\STATE{Compress the computed PGD solution: \\
{\small \verb|vi0 = pgdCompression(vi0,`tolModes',|}$\varepsilon^\star${\small \verb|});|}
}
\ENDIF
\ENSURE{Solution $\bvpgd_{i,0}$ stored in the form of a \texttt{separatedTensor} structure containing the vectors $\bV_{i,0}^m$ and $\bphi_{i,0}^m$ of the spatial and parametric modes.}
\end{algorithmic}
\end{algorithm}
\begin{algorithm}[!ht]
\caption{Encapsulated PGD solver for problem~\eqref{eq:greedyB}.}\label{alg:setupB}
\begin{algorithmic}[1]
\REQUIRE{Spatial ($\mat{K}_{\nu}^1, \mat{K}_{\nu}^2, \mathbf{f}_\Lambda^{q, 1}, \mathbf{f}_\Lambda^{q, 2}$) and parametric ($\bxi_{\nu}^1, \bxi_{\nu}^2, \bups_{\Lambda}^1, \bups_{\Lambda}^2$) modes of the linear system, tolerances $\varepsilon$ for the PGD enrichment and $\varepsilon^\star$ for the PGD compression.}
\STATE{Initialization of the separated tensor for the left-hand side of problem: \\
{\small \verb|Kij = separatedTensor;|}
}
\STATE{Setup of the spatial modes: \\
{\small \verb|Kij.sectionalData{1,1} =|} $\mat{K}_{\nu}^1${\small \verb|;|} \\
{\small \verb|Kij.sectionalData{1,2} =|} $\mat{K}_{\nu}^2${\small \verb|;|}
}
\STATE{Setup of the parametric modes depending on the physical parameter $\mu$: \\
{\small \verb|Kij.sectionalData{2,1} =|} $\bxi_{\nu}^1${\small \verb|;|} \\
{\small \verb|Kij.sectionalData{2,2} =|} $\bxi_{\nu}^2${\small \verb|;|} \\
}
\FOR{{\small \texttt{q = 1:}}card$(\mathcal{N}_i^j)$}
\STATE{Setup of the parametric modes depending on the active boundary parameters $\bLambda_i^j$: \\
{\small \verb|Kij.sectionalData{2+q,1} =|} $\bups_{\Lambda}^1${\small \verb|;|} \\
{\small \verb|Kij.sectionalData{2+q,2} =|} $\bups_{\Lambda}^1${\small \verb|;|} \\
}
\ENDFOR
\STATE{Initialization of the separated tensor for the right-hand side of the problem: \\
{\small \verb|Fij = separatedTensor;|}
}
\FOR{{\small \texttt{q = 1:}}card$(\mathcal{N}_i^j)$}
\STATE{Setup of the spatial modes: \\
{\small \verb|Fij.sectionalData{1,2*q-1} =|} $\mathbf{f}_\Lambda^{q, 1}${\small \verb|;|} \\
{\small \verb|Fij.sectionalData{1,2*q}   =|} $\mathbf{f}_\Lambda^{q, 2}${\small \verb|;|}
}
\STATE{Setup of the parametric modes depending on the physical parameter $\mu$: \\
{\small \verb|Fij.sectionalData{2,2*q-1} =|} $\bxi_{\nu}^1${\small \verb|;|} \\
{\small \verb|Fij.sectionalData{2,2*q}   =|} $\bxi_{\nu}^2${\small \verb|;|}
}
\FOR{{\small \texttt{r = 1:}}card$(\mathcal{N}_i^j)$}
\STATE{Setup of the parametric modes depending on the active boundary parameters $\bLambda_i^j$:}
\IF{{\small \texttt{r == q}}}
\STATE{
{\small \verb|Fij.sectionalData{2+r,2*q-1} =|} $\bups_{\Lambda}^2${\small \verb|;|} \\
{\small \verb|Fij.sectionalData{2+r,2*q}   =|} $\bups_{\Lambda}^2${\small \verb|;|}
}
\ELSE
\STATE{
{\small \verb|Fij.sectionalData{2+r,2*q-1} =|} $\bups_{\Lambda}^1${\small \verb|;|} \\
{\small \verb|Fij.sectionalData{2+r,2*q}   =|} $\bups_{\Lambda}^1${\small \verb|;|}
}
\ENDIF
\ENDFOR
\ENDFOR
\STATE{Solution of the problem via the alternating directions method: \\
{\small \verb|vij = pgdLinearSolve(Kij,Fij,`tolModes',|}$\varepsilon${\small \verb|});|}
}
\IF{PGD compression is active}
\STATE{Compress the computed PGD solution: \\
{\small \verb|vij = pgdCompression(vij,`tolModes',|}$\varepsilon^\star${\small \verb|});|}
}
\ENDIF
\ENSURE{Solution $\bvpgd_{i,j}$ stored in the form of a \texttt{separatedTensor} structure containing the vectors $\bV_{i,j}^m$, $\bphi_{i,j}^m$ and $\bpsi_{i,j}^m$ of the spatial and parametric modes.}
\end{algorithmic}
\end{algorithm}

\clearpage

\section{SUPG formulation on a reference domain}
\label{append:Graetz}

The problem in Sect.~\ref{sec:testRozza} features a convection-diffusion equation in a parametric domain. To construct the PGD surrogate model, a stabilized SUPG formulation is employed and the subproblem in $\Omega_2(\bmu) = [1, 1+\mu_2] \times [0,1]$ is rewritten on a parameter-independent subdomain $\hOmega_2 = [0,1] \times [0,1]$ via a reference domain configuration.

First, the variational form of problems~\eqref{eq:pbApgd} and~\eqref{eq:pbBpgd} with SUPG stabilization and $\mathbb{Q}_1$ Lagrangian finite elements are recalled. More precisely, the bilinear form is given by
\begin{equation}\label{eq:bilinearSUPG}
\begin{aligned}
\Apgd(v, \deV; \bmu) 
=& \frac{1}{\mu_1} \int_{\Omega_2(\bmu)}   \nabla v\cdot \nabla \deV\,d\bx 
+ \int_{\Omega_2(\bmu)}  \balpha \cdot \nabla v \, \deV\,d\bx 
\\
&+ \sum_{T \in \Omega_2(\bmu)} \int_{T} \tSUPG (\balpha \cdot \nabla v) (\balpha \cdot \nabla \deV)\,d\bx 
\, ,
\end{aligned}
\end{equation}
while the linear forms for problems~\eqref{eq:pbApgd} and~\eqref{eq:pbBpgd} respectively account for the parameter-independent Dirichlet boundary condition on $\Gamma^{D,2}_2(\bmu)$ and for the Dirichlet condition associated with the active boundary parameters on the interface.

Following~\cite{Huerta-GCCDH-13} and exploiting the horizontal direction of the studied convection field, the stabilization parameter is defined as
\begin{equation}\label{eq:tauSPUG}
    \tSUPG = h_{x,1} \left(1 + \frac{9}{\mathrm{Pe}^2}\right)^{-\tfrac{1}{4 |\alpha_1|}}.
\end{equation}
Of course, the stabilization coefficient $\tSUPG$ depends upon the parameters. A detailed discussion on the choice of such stabilization in the context of PGD-ROM is available in~\cite{Huerta-GCCDH-13}. Nonetheless, in the present study, it was observed that the value of $\tSUPG$ computed according to the definition~\eqref{eq:tauSPUG} does not significantly vary with the parameters. Hence, 
for all the computations, a space-dependent stabilization coefficient is obtained by setting in~\eqref{eq:tauSPUG} the value of the P\'eclet number associated with $\mu_1 = 2 \times 10^4$.

\begin{figure}[h!]
    \centering
    \begin{tikzpicture}
        \draw (0, 0) node[anchor = north east]{$(0,0)$} -- (0.1, 0) node[anchor = north west]{$(\bar{h},0)$} -- (3, 0) node[anchor = north]{$(1, 0)$};
        \draw (0, 0) -- (0, 3) node[anchor = south east]{$(0, 1)$};
        \draw (0, 3) -- (0.1, 3) node[anchor = south west]{$(\bar{h},1)$} -- (3, 3) node[anchor = south]{$(1, 1)$};
        \draw (3, 0) -- (3, 3);
        \draw[dotted] (0.2, 0) -- (0.2, 3);

        \node at (1.5, 1.5) {$\hOmega_2$};

        \draw[->, thick] (3.2, 1.2) -- (5.8, 1.2);
        \node at (4.5, 1.5) {$\Map$};

        \draw (6, 0) node[anchor = north east]{$(1,0)$} -- (6.1, 0) node[anchor = north west]{$(1+\bar{h},0)$} -- (11, 0) node[anchor = north]{$(1+\mu_2, 0)$};
        \draw (6, 0) -- (6, 3) node[anchor = south east]{$(1, 1)$};
        \draw (6, 3) -- (6.1, 3) node[anchor = south west]{$(1+\bar{h},1)$} -- (11, 3) node[anchor = south]{$(1+\mu_2, 1)$};
        \draw (11, 0) -- (11, 3);
        \draw[dotted] (6.2, 0) -- (6.2, 3);
        
        \node at (8.5, 1.5) {$\Omega_2$};
    \end{tikzpicture}
    \caption{Definition of the parametric mapping $\Map$.}
    \label{fig:mapGraetz}
\end{figure}
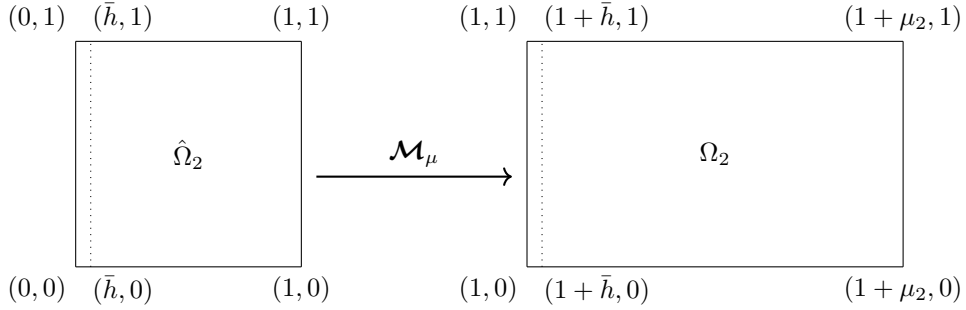

In order to construct a PGD surrogate model starting from the SUPG formulation presented above, all terms in the bilinear and linear forms need to be appropriately rewritten in a parameter-independent domain.
To this end, the parametric mapping
\begin{equation}\label{eq:mapping}
\begin{aligned}
\Map : \
& \hOmega_2 \times \mathcal{I}^2 \rightarrow \Omega_2(\bmu) \\[0.5em]
& (\hx,\hy,\mu_2) \mapsto (x,y)
\end{aligned}
\end{equation}
between the reference and the physical subdomains is defined as
\begin{equation}\label{eq:transformation}
x = 
\begin{cases}
    1+ \hx \quad & \text{for } \hx \leq \bar{h}\,,\\[0.5em]
    \displaystyle\frac{1 - \bar{h} \hx}{1-\bar{h}} + \mu_2 \frac{\hx - \bar{h}}{1-\bar{h}}\quad &\text{for } \hx > \bar{h} \, ,
\end{cases}
\qquad
y = \hy \, ,
\end{equation}
with $\bar{h} = 5 \times 10^{-2}$. A sketch of the transformation is displayed in Figure~\ref{fig:mapGraetz}.

Hence, following~\cite{Ammar-AHCCL-14}, the integrals in~\eqref{eq:bilinearSUPG} are mapped to the reference subdomain by inverting the transformation~\eqref{eq:mapping}. Let $\Jaco$ denote the Jacobian of the mapping, $\detJ$ its determinant and $\adjJ = \detJ \Jaco^{-1}$ its adjoint.
The resulting SUPG bilinear form on the reference subdomain is given by
\begin{equation}\label{eq:bilinearSUPGfixed}
\begin{aligned}
\Apgd(v, \deV; \bmu) 
=& \frac{1}{\mu_1} \int_{\hOmega_2}   \nabla v\cdot \left(\frac{\adjJt \adjJ}{\detJ} \nabla \deV \right) d\bhx 
+ \int_{\hOmega_2}  \balpha \cdot \left( \adjJ \nabla v \right) \, \deV\,d\bhx 
\\
&+ \sum_{\hT \in \hOmega_2} \int_{\hT} \tSUPG (\balpha \cdot \adjJ \nabla v) (\balpha \cdot \Jaco^{-1}\nabla \deV)\,d\bhx 
\, ,
\end{aligned}
\end{equation}
where all integrals are computed on the parameter-independent subdomain $\hOmega_2$ and the dependence on $\mu_2$ is encapsulated in the mapping.
Introducing the definition of $\Jaco$, $\detJ$ and $\adjJ$ for the transformation~\eqref{eq:transformation} into~\eqref{eq:bilinearSUPGfixed}, the resulting bilinear form features an affine dependence on $\bmu$ (see, e.g.,~\cite{Rozza:14}), yielding
standard finite element matrices with non-constant parameters. Hence, the strategy described in Appendix~\ref{append:Setup} can be straightforwardly applied to equation~\eqref{eq:bilinearSUPGfixed} to construct the local surrogate model using the encapsulated PGD.
For alternative strategies to construct PGD surrogate models of geometrically parametrized problems featuring more general transformations, the interested reader is referred to~\cite{Zlotnik:2015:IJNME,Sevilla-SZH-20,Sevilla-SBGH-20}.

	
\end{document}